\documentclass[hidelinks,onefignum,onetabnum]{siamart220329}

\usepackage{hyperref}
\hypersetup{pdfauthor={Name}}
\hypersetup{colorlinks=true}

\usepackage{amsfonts}
\usepackage{amsmath}
\usepackage{subcaption}
\usepackage{tikz}
\usepackage{graphicx}

\usepackage{algorithm}
\usepackage{algpseudocode}

\usepackage{bm}         %
\usepackage{siunitx}    %
\usepackage{booktabs} %
\usepackage{multirow}   %
\usepackage{placeins} %
\usepackage{soul}
\usepackage{mathtools} %

\usepackage[most]{tcolorbox}

\ifpdf
  \DeclareGraphicsExtensions{.eps,.pdf,.png,.jpg}
\else
  \DeclareGraphicsExtensions{.eps}
\fi

\usepackage{nameref}

\newtcbox{\highlight}[2]{enhanced, box align=base, nobeforeafter, colback=#1,
colframe=#2, size=small, left=0pt, right=0pt, boxsep=1.6pt, boxrule=1.2pt}
\newcommand{\dashedredframe}[1]{\tikz[baseline=(char.base)]{\node[anchor=text,
rectangle, draw=black!60!red, dashed, dash pattern= on 0.5pt off 1.0pt, thin, inner sep=2pt]     (char) {#1};}}
\newcommand{\blueframe}[1]{\tikz[baseline=(char.base)]{\node[anchor=text,
rectangle, draw=black!60!blue, inner sep=2pt] (char) {#1};}}

\newcommand{\githuburl}{https://github.com/Alexey-Voronin/Monolithic_AMG_For_Stokes}
\newcommand{\githuburltext}{https://github.com/Alexey-Voronin/Monolithic\_AMG\_For\_Stokes}
\renewcommand{\vec}[1]{\bm{#1}}
\newcommand{\qtwo}{$\mathbb{Q}_2$}

\newcommand{\pone}{\ensuremath{\mathbb{P}_1}}
\newcommand{\ponedisc}{\ensuremath{\mathbb{P}_1^{disc}}}
\newcommand{\isoqtwo}{$\mathbb{Q}_1 \text{iso}\kern1pt\mathbb{Q}_2$}
\newcommand{\isoptwo}{$\mathbb{P}_1 \text{iso}\kern1pt\mathbb{P}_2$}
\newcommand{\vecisoptwo}{$\pmb{\mathbb{P}}_1 \text{iso}\kern1pt\pmb{\mathbb{P}}_2$}
\newcommand{\qtwoqone}{$\pmb{\mathbb{Q}}_2/\mathbb{Q}_1$}
\newcommand{\ptwopone}{$\pmb{\mathbb{P}}_2/\mathbb{P}_1$}
\newcommand{\ponepone}{$\pmb{\mathbb{P}}_1/\mathbb{P}_1$}
\newcommand{\isoqtwoqone}{$\pmb{\mathbb{Q}}_1 \text{iso}\kern1pt\pmb{ \mathbb{Q}}_2/             \mathbb{Q}_1$}
\newcommand{\isoptwopone}{$\pmb{\mathbb{P}}_1 \text{iso}\kern1pt\pmb{ \mathbb{P}}_2/             \mathbb{P}_1$}
\newcommand{\ptwoponedisc}{$\pmb{\mathbb{P}}_2/\mathbb{P}_1^{disc}$}
\newcommand{\pkpkmonedisc}{$\pmb{ \mathbb{P}}_k/\mathbb{P}_{k-1}^{disc}$}
\newcommand{\pkpkmone}{$\pmb{ \mathbb{P}}_k/\mathbb{P}_{k-1}$}
\newcommand{\pk}{$\pmb{ \mathbb{P}}_k$}
\newcommand{\pkmonedisc}{$\mathbb{P}_{k-1}^{disc}$}

\newcommand{\Mh}{HO-AMG}
\definecolor{colorA}{rgb}{0.12156862745098039, 0.4666666666666667, 0.7058823529411765}
\definecolor{colorB}{rgb}{1.0, 0.4980392156862745, 0.054901960784313725}
\definecolor{colorC}{rgb}{0.17254901960784313, 0.6274509803921569, 0.17254901960784313}
\definecolor{colorD}{rgb}{0.8392156862745098, 0.15294117647058825, 0.1568627450980392}
\definecolor{colorE}{rgb}{0.5803921568627451, 0.403921568627451,
0.7411764705882353}

\newcommand{\cMBTP}{{\color{colorE} BTP}} %
\newcommand{\cMh}{{\color{colorA} HO-AMG}} %
\newcommand{\cMphHO}{{\color{colorB} DC-skip1}} %
\newcommand{\cMphLO}{{\color{colorC} DC-skip0}} %
\newcommand{\cMphHLO}{{\color{colorD} DC-all}}

\newcommand{\cMphHOVanka}{{\color{colorB} DC-skip1(Vanka)}}
\newcommand{\cMphHODLSC}{{\color{colorB} DC-skip1(DLSC)}}
\newcommand{\cMphHLOVanka}{{\color{colorD} DC-all(Vanka)}}
\newcommand{\cMphHLOLSC}{{\color{colorD} DC-all(Vanka-DLSC)}}
\Crefname{algocf}{Algorithm}{Algorithms}
\newtheorem{remark}[theorem]{Remark}

\headers{Monolithic AMG for Stokes Equations
}{A. Voronin, S. MacLachlan, L.N. Olson, and R.S. Tuminaro}

\title{Monolithic Algebraic Multigrid Preconditioners for the Stokes 
 Equations\thanks{Submitted to the editors 07/07/2023.
\funding{The work of SPM was partially supported by an NSERC Discovery Grant. RT was
 supported by the U.S.~Department of Energy, Office of Science, Office of
 Advanced Scientific Computing Research, Applied Mathematics program. }}}

\author{Alexey Voronin\thanks{Department of Computer Science, University of
Illinois at Urbana-Champaign, Urbana IL, 61801 USA 
  (\email{voronin2@illinois.edu}, \url{http://alexey-voronin.github.io/},
  \email{lukeo@illinois.edu}, \url{http://lukeo.cs.illinois.edu/}
  ).}
\and Scott MacLachlan\thanks{Department of
 Mathematics and Statistics, Memorial University of Newfoundland, St.~John's NL, Canada 
  (\email{smaclachlan@mun.ca}, \url{https://www.math.mun.ca/\string~smaclachlan/}).
  }
\and Luke N. Olson \footnotemark[2]
\and Raymond S. Tuminaro\thanks{Sandia National Laboratories, Livermore, CA 94551 USA(\email{rstumin@sandia.gov}, \url{https://www.sandia.gov/ccr/staff/raymond-s-tuminaro/}).
  }
}

\usepackage{amsopn}

\algrenewcommand{\algorithmiccomment}[1]{\hfill{\scriptsize\{#1\}}}

\ifpdf
\hypersetup{
  pdftitle={Monolithic Algebraic Multigrid Preconditioners for the Stokes
 Equations},
  pdfauthor={A. Voronin, S. MacLachlan, L. N. Olson, and R. Tuminaro}
}
\fi

\begin{document}

\maketitle

\begin{abstract}
We investigate a novel monolithic algebraic multigrid (AMG)
preconditioner for the Taylor-Hood (\ptwopone{}) and Scott-Vogelius (\ptwoponedisc{}) discretizations of the Stokes equations. The algorithm is based on the use of the lower-order
    \isoptwopone{} operator within a defect-correction setting, in combination with
    AMG construction of
    interpolation operators for velocities and pressures. The 
    preconditioning framework is primarily algebraic, though the \isoptwopone{} 
    operator must be provided.
    We investigate two relaxation strategies in this setting.
    Specifically, a novel block factorization approach is devised for Vanka patch systems, which 
    significantly reduces storage requirements and computational overhead, and a Chebyshev
    adaptation of the LSC-DGS relaxation from~\cite{wang2013multigrid} is developed to improve parallelism.
The preconditioner demonstrates robust performance across a variety of 2D and 3D
Stokes problems, often matching or exceeding the effectiveness of an 
inexact block-triangular (or Uzawa) preconditioner,
especially in challenging scenarios such as elongated-domain problems. 
\end{abstract}



\section{Introduction}\label{sec:intro}

Differential equations and discrete problems with saddle-point structure arise
in many scientific and engineering applications. Common examples are the
incompressible Stokes and Navier-Stokes equations, where the incompressibility
constraint leads to saddle-point problems both at the continuum level and for
many common discretizations~\cite{elman2014finite,MBenzi_GHGolub_JLiesen_2005a}.
Discretizations of these systems using stable finite-element pairs lead to
$2\times2$ block structured linear systems, where the
lower-right diagonal block is a zero-matrix. The indefiniteness in
these matrices, along with the coupling between the two unknown types,
prevents the straightforward application of standard geometric or algebraic
multigrid methods.  Several successful preconditioners exist for these systems 
and are typically
built on approximate block-factorizations~\cite{elman2014finite, adler2017preconditioning, farrell2021reynolds}
or monolithic geometric multigrid~\cite{ABrandt_NDinar_1979a,
SPVanka_1986a, JLinden_etal_1989a, adler2017preconditioning, Braess,
schoberl2003schwarz}. While some algebraic multigrid algorithms have been
proposed~\cite{MWabro_2004a, MWabro_2006a, metsch2013algebraic,
AJanka_2008a, prokopenko2017algebraic}, general-purpose algebraic multigrid
methods that are successful for the numerical solution of a wide variety of
discretizations of the Stokes equations are recent (and still rare)~\cite{notay2016new,notay2017algebraic,bacq2022new,
bacq2021all}.  In this work, we propose a robust and
efficient algebraic multigrid method for the Stokes \isoptwopone{} discretization
and show how the method can be leveraged to develop effective
preconditioners for other discretizations, building on the geometric
multigrid methods in~\cite{voronin2021low}.

Monolithic geometric multigrid methods for these problems are well-known and
date back to some of the earliest multigrid approaches for
finite-difference
discretizations~\cite{ABrandt_1984a,SPVanka_1986a,JLinden_etal_1989a,ANiestegge_KWitsch_1990a}.
In these works and in later works focused on finite-element
discretizations~\cite{VJohn_LTobiska_2000a, john2002non, MLarin_AReusken_2008a,
emamilehrstuhl, BGmeiner_etal_2016a,
adler2017preconditioning,adler2021monolithic}, simultaneous geometric coarsening
of velocity and pressure fields leads to coarse-grid operators that naturally
retain stability on all grids of the hierarchy.  This coarsening in combination
with monolithic relaxation schemes, such as Vanka~\cite{SPVanka_1986a},
Braess-Sarazin~\cite{Braess}, or distributive~\cite{ABrandt_NDinar_1979a,wang2013multigrid}
relaxation, has led to effective multigrid convergence for a wide variety of
problems and discretizations. Despite the success and effectiveness of these methods,
they are not without their challenges, particularly when considering their extension to
algebraic multigrid approaches.

Two difficulties arise when attempting
to extend these monolithic methods to algebraic multigrid approaches.  First, it is
difficult to construct independent algebraic coarsenings of the velocity and
pressure fields that maintain algebraic properties of the fine-grid
discretization (such as the relative number of velocity and pressure degrees of
freedom) after Galerkin coarsening.  As a result, it is possible that a stable
discretization on the finest grid yields unstable discretizations on coarser
grids in the algebraic multigrid hierarchy.  Secondly, most inf-sup stable
finite element discretizations of the Stokes equations rely on higher-order
finite-element spaces (particularly for the discretized velocity), yielding
discretization matrices that are far from the M-matrices for which classical AMG
heuristics are known to work best~\cite{heys2005algebraic, olson2007algebraic}.
Consequently, even if stability problems are avoided, convergence of AMG-based
solvers can deteriorate as the order of the discretization increases.

Many existing AMG approaches for the Stokes equations attempt to
bypass coarsening issues by constraining
the coarsening within AMG to loosely mimic properties of geometric multigrid.
Specifically, it is possible to roughly maintain the ratio between the number of
velocity degrees of freedom and the number of pressure degrees of freedom throughout the multigrid hierarchy
by leveraging either geometric physical degree of freedom locations or element
information~\cite{MWabro_2004a,MWabro_2006a,AJanka_2008a,metsch2013algebraic,prokopenko2017algebraic}.
While successful, this introduces additional heuristics that depend on both the
fine-level discretization and the geometry of the mesh in the AMG process.
Change-of-variable transformations have been used successfully to construct
monolithic AMG preconditioners for saddle-point-type systems such as the Stokes
equations~\cite{webster2013stability,notay2016new,notay2017algebraic,bacq2022new,
bacq2021all}. These transformations convert the system into a similar one with a
scalar elliptic operator in the bottom-right block, making it more suitable for
classical AMG with unknown-based coarsening and point-wise relaxation. The
method has been theoretically analyzed, showing a uniform bound on the spectral
radius of the iteration matrix using a single step of damped Jacobi
relaxation~\cite{notay2016new}. However, when applying these methods to Oseen
problems using a two-level approach, it has been demonstrated that the use of
more effective relaxation approaches, such as Gauss-Seidel, is necessary to achieve
fast convergence~\cite{bacq2022new}. While this approach has been
extended to multilevel contexts, as demonstrated in~\cite{bacq2021all}, ensuring robust and cost-effective convergence requires
two pieces: first, the efficient implementation of relaxation techniques like Gauss-Seidel
and successive over-relaxation (SOR), which are difficult to parallelize; and
second, effective algorithms for sparsification of coarse-grid matrices.

In this paper, we propose a new algebraic multigrid framework for the solution
of the discretized Stokes equations, based on preconditioning higher-order
discretizations by smoothed-aggregation AMG applied to the low-order
\isoptwopone{} discretization. We employ both Vanka-style relaxation and a generalization of the least-squares commutator distributive Gauss-Seidel (LSC-DGS) method from~\cite{wang2013multigrid}.
These are applied to the original high-order discretization and each level
within the algebraic multigrid (AMG) hierarchy, allowing us to compare their
relative performances.  The use of a low-order preconditioner (in a classical
multigrid \textit{defect correction} approach) allows us to avoid issues that
arise when applying AMG directly to higher-order discretizations, leading to
robust convergence with low
computational complexities.  In~\cite{voronin2021low}, we considered a similar
strategy within geometric multigrid for the \qtwoqone{} discretization, using
defect correction based on the \isoqtwoqone{} discretization.  Here, we extend
that work to both structured and unstructured triangular and tetrahedral meshes,
using both Taylor-Hood, \ptwopone{}, and Scott-Vogelius, \ptwoponedisc{},
discretizations on the finest grids.  To achieve scalable results, particularly
on unstructured grids, we examine best practices for algebraic multigrid on the
\isoptwopone{} discretization, including the preservation of stability
on coarse levels of the hierarchy.  We note that this also provides insight into
achieving scalable results applying AMG directly to the \ptwopone{}
discretization, and we compare this to our defect correction approach.

The remainder of this paper is structured as follows.
In~\cref{sec:disc}, we introduce the Stokes equations and the accompanying
discretizations. \Cref{sec:mg_setup} describes the monolithic AMG
framework and the construction of the low-order preconditioner for high-order
discretizations of the Stokes equations. \Cref{sec:results} presents numerical
results, demonstrating the robustness of the preconditioner for a set of
test problems, as well as a comparison of the computational cost between the
various AMG-based preconditioners. In~\cref{sec:application_problems},
we present numerical results obtained by applying the most effective
preconditioners identified in the previous section to more complex application
problems. \Cref{sec:conclusion} provides concluding
remarks, with future research directions outlined in \cref{sec:future_work}.

\section{Discretization of the Stokes Equations}\label{sec:disc}

In this paper, we consider the steady-state Stokes problem on a bounded
Lipschitz domain $\Omega \subset \mathbb{R}^d$, for $d=2$ or $3$, with boundary
$\partial\Omega=\Gamma_{\text{N}} \cup \Gamma_{\text{D}}$, given by
\begin{subequations}\label{eq:stokes-eq}
  \begin{alignat}{2}
    -\nabla^2 \vec{u} +\nabla p &=\vec{f}              &\quad& \text{in $\Omega$} \label{eq:stokes-eq1} \\
    -\nabla \cdot \vec{u}       &= 0                   &\quad& \text{in $\Omega$} \label{eq:stokes-eq2} \\
                        \vec{u} &= \vec{g}_{\text{D}}  &\quad& \text{on $\Gamma_{\text{D}}$} \label{eq:stokes-eq3} \\
        \frac{\partial\vec{u}}{\partial \vec{n}} - \vec{n} p
                                &= \vec{g}_{\text{N}}  &\quad& \text{on $\Gamma_{\text{N}}$}. \label{eq:stokes-eq4}
\end{alignat}
\end{subequations}
Here, $\vec{u}$ is the velocity of a viscous fluid, $p$ is the pressure,
$\vec{f}$ is a forcing term, $\vec{n}$ is the outward pointing normal to
$\partial\Omega$,
$\vec{g}_{\text{D}}$ and $\vec{g}_{\text{N}}$ are given boundary data, and
$\Gamma_{\text{D}}$ and $\Gamma_{\text{N}}$ are disjoint Dirichlet and Neumann
segments of the boundary.  In a variational formulation, the natural function
space for the velocity is $\vec{\mathcal{H}}_{\vec{g}_{\text{D}}}^1(\Omega)$,
the subset of $\vec{\mathcal{H}}^1(\Omega)$ constrained to match the essential
boundary condition in~\eqref{eq:stokes-eq3}.  If the velocity is specified
everywhere along the
boundary, which requires $\int_{\partial\Omega} \vec{g}_{\text{D}}\cdot\vec{n}=0$,
then a suitable function space for the pressure is
the space of zero-mean functions in $L^2(\Omega)$~---~i.e., $\mathcal{Q}=L^2_0(\Omega)$;
this avoids a pressure solution that is unique only up to a constant.
In the case where $|\Gamma_{\text{N}}| > 0$,
this additional constraint is not needed, and the pressure is assumed to be in
$\mathcal{Q}=L^2(\Omega)$~\cite{elman2014finite}.

For both Taylor-Hood (TH) and Scott-Vogelius (SV) element pairs, define
finite-dimensional spaces $(\vec{\mathcal{V}}_h^{\text{D}},
\mathcal{Q}_h)\subset (\vec{\mathcal{H}}_{\vec{g}_{\text{D}}}^1(\Omega),
\mathcal{Q})$, where $\vec{\mathcal{V}}_h^{\text{D}}$ strongly satisfies the
Dirichlet boundary conditions on $\Gamma_{\text{D}}$ specified
by~\cref{eq:stokes-eq3} (at least up to interpolation error).  The resulting
weak formulation of~\cref{eq:stokes-eq} is to find $\vec{u} \in
\vec{\mathcal{V}}_h^{\text{D}}$ and $p \in \mathcal{Q}_h$ such that
\begin{subequations}\label{eq:stokes-disc}
  \begin{alignat}{2}
    & a(\vec{u}, \vec{v})+b(p, \vec{v}) &&= F(\vec{v})\\
    & b(q, \vec{u})                     &&= 0,
  \end{alignat}
\end{subequations}
for all $\vec{v} \in \vec{\mathcal{V}}_h^0$ and $q \in \mathcal{Q}_h$, where
$\vec{\mathcal{V}}_h^0$ is the same finite-element space as
$\vec{\mathcal{V}}_h^{\text{D}}$, but with zero Dirichlet boundary conditions on
$\Gamma_{\text{D}}$.  Here, $a(\cdot,\cdot)$ and $b(\cdot,\cdot)$ are bilinear
forms, and $F(\cdot)$ is
    a linear form given by
\begin{equation*}
        a(\vec{u}, \vec{v}) = \int_{\Omega} \nabla \vec{u} : \nabla \vec{v},\quad
        b(p, \vec{v}) = - \int_{\Omega} p \nabla \cdot \vec{v},\quad\text{and}\quad
        F(\vec{v}) = \int_{\Omega} \vec{f} \cdot \vec{v}+\int_{\Gamma_{\text{N}}} \vec{g}_{\text{N}} \cdot \vec{v}.
\end{equation*}
From this point forward, we overload the notation and use $\vec{u}$ and
$p$ to denote the discrete velocity and pressure unknowns in a finite-element
discretization of~\cref{eq:stokes-disc}.  Given an inf-sup stable choice of
finite-dimensional spaces $(\vec{\mathcal{V}}_h^{\text{D}}, \mathcal{Q}_h)$, we obtain a saddle-point system of the form
\begin{equation}\label{eq:saddle}
    Kx
    =
    \begin{bmatrix}
        A  & B^T \\
        B  &  0
    \end{bmatrix}
    \begin{bmatrix}
        \vec{u}\\
        p
    \end{bmatrix}
    =
    \begin{bmatrix}
        \vec{f}\\
        0
    \end{bmatrix}
    =b,
\end{equation}
where matrix $A$ corresponds to the discrete vector Laplacian on
$\mathcal{V}_h^{\text{D}}$, and $B$ represents the negative of the discrete
divergence operator mapping $\mathcal{V}_h^{\text{D}}$ into $\mathcal{Q}_h$.
In two dimensions, the dimensions of $A$ and $B$ are given by $A \in
\mathbb{R}^{(n_{x}+n_{y})\times (n_{x}+n_{y})}$ and $B \in \mathbb{R}^{n_{p}\times
(n_{x}+n_{y})}$, where $n_x$ and $n_y$ represent the number of velocity nodes
for each component of the vector $\vec{u}$ (which are equal for the boundary
conditions considered here) and $n_p$ represents the number of pressure nodes.

\subsection{Finite Element Discretizations and Meshes}

For the remainder of the paper, we focus on three different but related
discretizations to define finite-dimensional spaces $(\vec{\mathcal{V}}_h^{\text{D}},
\mathcal{Q}_h)$ on meshes of $d$-dimensional simplices. The first is the standard
Taylor-Hood (TH) discretization, \pkpkmone{}, with continuous piecewise
polynomials of degree $k$ for the velocity space $\vec{\mathcal{V}}_h^{\text{D}}$ and
continuous piecewise polynomials of degree $k-1$ for the pressure space
$\mathcal{Q}_h$. The TH element pair is inf-sup stable for $k\geq2$ on any
triangular (tetrahedral) mesh $\Omega_h$ of
domain $\Omega$. Here, we focus on the case $k=2$, leading to the lowest-order pair of
TH elements, corresponding to piecewise quadratic
velocity components and piecewise linear elements for the pressure.

\begin{figure}[!ht]
  \centering
  \begin{subfigure}{0.2\textwidth}
    \includegraphics[width=\textwidth]{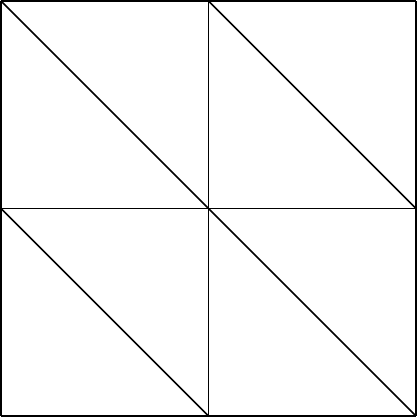}
    \caption{$\Omega_h$}\label{fig:mesh_a}
  \end{subfigure}
  \hfill
  \begin{subfigure}{0.2\textwidth}
    \includegraphics[width=\textwidth]{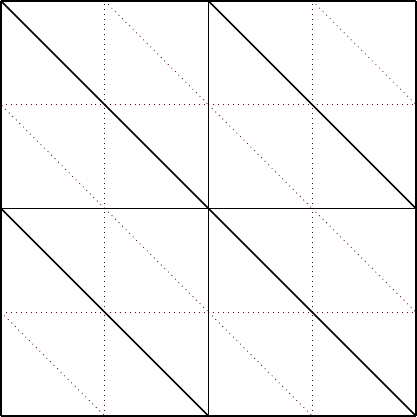}
    \caption{$R_\mathcal{U}(\Omega_h)$}\label{fig:mesh_b}
  \end{subfigure}
  \hfill
  \begin{subfigure}{0.2\textwidth}
    \includegraphics[width=\textwidth]{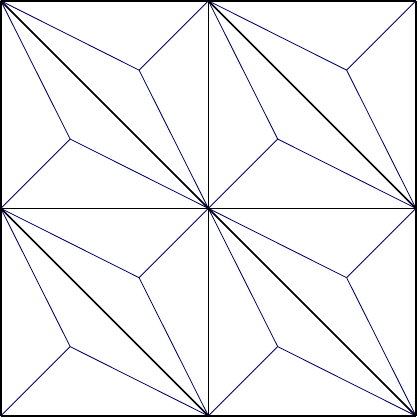}
    \caption{$R_\mathcal{B}(\Omega_h)$}\label{fig:mesh_c}
  \end{subfigure}
  \hfill
  \begin{subfigure}{0.2\textwidth}
    \includegraphics[width=\textwidth]{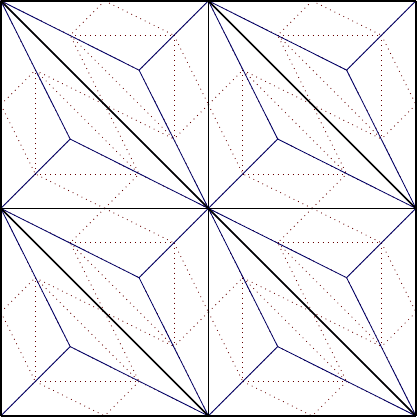}
    \caption{$R_\mathcal{U}(R_\mathcal{B}(\Omega_h))$}\label{fig:mesh_d}
  \end{subfigure}
  \captionsetup{skip=1pt}
  \caption{Mesh refinements.
    (a) structured 2D triangular mesh;
    (b) the \protect\dashedredframe{quadrisection refinement} of (a);
    (c) barycentric mesh (Alfeld split), obtained via
    \protect\blueframe{barycentric refinement} of (a); and
    (d) \protect\dashedredframe{quadrisection refinement} of (c).}\label{fig:meshes}
\end{figure}

The second discretization that we consider is the Scott-Vogelius (SV) discretization,
\pkpkmonedisc{}, comprised of continuous piecewise polynomials of degree $k$ for
the velocity and discontinuous piecewise polynomials of degree $k-1$ for the
pressure. The inf-sup stability of the SV discretization depends on the
polynomial degree ($k$) of the function spaces and properties of the mesh used
to assemble the system. In~\cite{qin1994convergence,zhang2005new}, the
\pkpkmonedisc{} discretization is shown to be stable for $k \geq d$ on any mesh
that results from a single step of barycentric refinement of any triangular
(tetrahedral) mesh.  Such a refinement is depicted in~\cref{fig:mesh_c}, where
$R_\mathcal{B}(\Omega_h)$ denotes the barycentric refinement of mesh
$\Omega_h$. In~\cite{zhang2008p1,zhang2011quadratic}, the polynomial order
required for inf-sup stability of the SV discretization is lowered to $k \geq
d-1$ when the mesh is given by a Powell-Sabin split of a triangular or
tetrahedral mesh, $\Omega_h$, but we do not consider such meshes. An
important advantage of the SV discretization over the TH discretization is that
the choice of the discontinuous discrete pressure space satisfies the
relationship that the divergence of any velocity in \pk{} is exactly represented
in the pressure space, \pkmonedisc{}; therefore, enforcing the divergence
constraint weakly leads to velocities that are point-wise divergence free as
well. This is not the case for TH elements, although an additional grad-div
stabilization term can be used to achieve stronger point-wise mass
conservation~\cite{john2017divergence,case2011connection,linke2011convergence}.
This work focuses on the \ptwoponedisc{} discretization on barycentrically
refined meshes in two dimensions. %

The third discretization is the low-order \isoptwopone{}
discretization. It is well known that a mixed discretization using piecewise linear
finite elements for both velocities and pressures, such as \ponepone{},
does not satisfy the necessary
inf-sup stability criterion for well-posedness.  The \isoptwopone{} gets around
this instability by considering first-order discretizations of both velocities
and pressures, but on different meshes.  Specifically, the
$\pmb{\mathbb{P}}_1$ velocity space is discretized on a quadrisection
refinement, $R_\mathcal{U}(\Omega_h)$, of the mesh, $\Omega_h$, used for the
$\mathbb{P}_1$ discretization of the pressure.  The \isoptwopone{}
discretization is inf-sup stable~\cite{ern2004book} but is not
commonly used due to its relatively low accuracy and computational efficiency.
Here, we make use of this discretization not to directly approximate solutions
to the Stokes equations, but as an auxiliary operator with which we
build preconditioners for the discretizations of interest.

\Cref{fig:meshes} illustrates meshes used in the construction of the \ptwopone{},
\ptwoponedisc{}, and \isoptwopone{} systems.  In all cases, we consider a given
mesh, $\Omega_h$, as shown in~\cref{fig:mesh_a}.  We denote the quadrisection
refinement of $\Omega_h$, by $R_\mathcal{U}(\Omega_h)$, as shown
in~\cref{fig:mesh_b}.  Similarly, \cref{fig:mesh_c} denotes the barycentric
refinement of $\Omega_h$, $R_\mathcal{B}(\Omega_h)$.  For the SV discretization
on mesh $R_\mathcal{B}(\Omega_h)$, we will make use of a quadrisection refinement of
$R_\mathcal{B}(\Omega_h)$, denoted $R_\mathcal{U}(R_\mathcal{B}(\Omega_h))$,
shown in~\cref{fig:mesh_d}.

\section{Monolithic Multigrid}\label{sec:mg_setup}

In this section, we introduce our multilevel algorithm for solving
the Stokes system, $Kx=b$, in~\cref{eq:saddle}.  The algorithm comprises several
key components that we detail in subsequent sections, namely:
\begin{itemize}
  \item the use of the low-order \isoptwopone\ discretization to precondition
    the system;
  \item the construction of algebraic multigrid hierarchies for separate
    portions of the system, including the vector Laplacian on the velocity
    variables (represented by $A$) and the \textit{div-grad}-like system for
        pressure;
  \item a careful algebraic coarsening of the $(\vec{u},p)$
    degrees-of-freedom (DoFs) to mimic behavior observed in efficient geometric
    multigrid schemes for Stokes;
  \item an algebraic form of Vanka relaxation for saddle-point systems, with a
    tuned relaxation parameter;
  \item a block LU factorization technique for the solution of Vanka patches that greatly reduces storage and computational cost; and
  \item a variant of the least-squares commutator distributive Gauss-Seidel (LSC-DGS) relaxation as an alternative to Vanka relaxation.
\end{itemize}
In total, the algorithm derived in this section (and presented in detail
in~\cref{alg:stokes_mg_setup}) is algebraic in fashion, requiring only a
description of the original (fine) DoFs and a low-order interpolation
operator based on the finite-element basis (on the fine level).

\subsection{Low-Order Preconditioning}

In~\cite{voronin2021low}, we demonstrate that the low-order \isoqtwoqone{}
finite-element discretization can be leveraged to construct efficient and robust
geometric multigrid (GMG) preconditioners for the higher-order \qtwoqone{}
discretization of the Stokes equations.  In that work, a multilevel
defect-correction approach was proposed that we follow here.  There are two aims
in the present work.  First, to extend the approach of~\cite{voronin2021low}
from two-dimensional quadrilateral meshes to triangular and tetrahedral meshes,
and to consider more discretizations.
Secondly, we aim to construct the
defect-correction preconditioner in an \textit{algebraic} setting without
relying on geometric information about the mesh to construct interpolation (or
other) operators.     While the method is described as an
AMG solver, we emphasize that the procedure is best used
as a preconditioner within a Krylov method.

In the algorithms that follow, $K_0$ and $K_1$ are the Stokes operators assembled using a higher-order
discretization (either \ptwopone{} or \ptwoponedisc{}) and the lower-order
\isoptwopone{} discretization, respectively, where the pressure for the
\isoptwopone{} discretization is defined on the same grid as the higher-order
discretization, and the \vecisoptwo{} velocity is defined on the quadrisection refinement
of this mesh.  As the spaces associated with $K_0$ and $K_1$ are different,
transfer operators are needed to map from the low-order mesh to the high-order
mesh, denoted by $P_0$, and vice-versa, denoted by $R_0$.  Further, let $G_1$
denote the error-propagation operator for a stationary iterative solver for
$K_1$. Extending~\cite{voronin2021low}, we consider a
defect-correction algorithm defined by the two-level iteration scheme with the
error-propagation operator
\begin{equation}\label{eq:TG-operator-full}
G_{dc} =\left(I - \omega_0 M_0^{-1} K_0\right)^{\nu_2}
      \left(I - P_0(I-G_1^{\gamma}) K_1^{-1} \eta R_0 K_0\right)
      \left(I - \omega_0 M_0^{-1} K_0\right)^{\nu_1}.
\end{equation}
Here, $\nu_1$ relaxation sweeps are first applied to the $K_0$ operator where
the specific relaxation scheme is defined by the choice of $M_0^{-1}$. This is
followed by a coarse/low-order correction given by the middle expression,
followed by post-relaxation on $K_0$. The low-order  correction includes
first a high-order residual calculation (the $K_0$ term), which is restricted using $\eta
R_0$. Then, $\gamma$ iterations of the $G_1$ solver are used to produce a
low-order approximate solution that is interpolated via $P_0$ back the fine
space. We note that if $\gamma = \infty$ and $G_1$ is a convergent iteration,
then the middle term reduces to $\left(I - P_0 K_1^{-1} \eta R_0 K_0\right)$,
which corresponds to solving the coarse
problem exactly.
As we demonstrate in \cref{sec:results}, an important (and non-standard) feature
of this cycle is the residual weighting operator, $\eta$.  Here, $\eta$ is a
$2\times 2$ block-diagonal matrix that imposes separate scaling of residual
equations corresponding to the momentum, \cref{eq:stokes-eq1}, and the
continuity, \cref{eq:stokes-eq2}, equations, by scalar weights $\eta_u$ and
$\eta_p$, and is
represented by the following operator:
\begin{equation}\label{eq:eta}
    \eta
    =
    \begin{bmatrix}
       \eta_u I  &  \\ %
               &  \eta_p I %
    \end{bmatrix}.
\end{equation}
Depending on the $K_0$ discretization, we find two useful
regimes for the scaling parameters, $\eta_u$ and $\eta_p$.  When $K_0$ comes
from the \ptwopone{} discretization, it is productive to use
$\eta_p\approx\eta_u \leq 1$, damping the defect correction from the $K_1$
discretization. In contrast, when $K_0$ comes from the
\ptwoponedisc{} discretization, we find $\eta_p \gg \eta_u \approx
1$, or continuity equation residual overweighting, is the most effective.  We
note that residual overweighting has been used to accelerate the
convergence of non-Galerkin defect-correction multigrid methods by minimizing
the discrepancy between levels in other
contexts~\cite{brandt1993accelerated,chen2015multigrid}.  We also note that, in the
work below, we do not introduce additional residual weighting between levels in
the multigrid cycle for $K_1$, as preliminary experiments found that this did not lead
to further improvements in convergence.

We construct the grid transfer operator, $P_0$, between the higher-order ($K_0$) and the low-order
($K_1$) discretization as a block-diagonal matrix
\begin{equation}\label{eq:interp}
    P_0
    =
    \begin{bmatrix}
        P^{\vec{u}}_0  &  \\
                     &  P^p_0
    \end{bmatrix},
\end{equation}
where $P^{\vec{u}}_0$ and $P^p_0$ correspond to velocity and pressure field
interpolation operators. For the $K_0$ and $K_1$ system
discretized with \ptwopone{} and \isoptwopone{} elements, respectively, the
pressure and velocity unknowns are co-located; therefore, it is plausible to use
injection, $P^{\vec{u}}_0=I$ and $P^{p}_0=I$, as in~\cite{voronin2021low}.
Formally, it is natural to define restriction using 
element projection (described shortly) %
to map between two different spaces.  However, given the
co-located nature of the degrees of freedom, it is convenient to simply use injection to define restriction
for \ptwopone{} and \isoptwopone{} discretizations, and so $R_0 = P_0^T$.
For the \ptwoponedisc{} operator, $K_0$, it is still convenient to use injection
for the velocity degrees of freedom, as the \vecisoptwo{} velocity degrees of
    freedom in $K_1$ are again co-located with those of $K_0$, again leading to
    $P^{\vec{u}}_0=R^{\vec{u}}_0=I$.  For pressure degrees of freedom,
    however, more complicated mappings are needed, as $K_0$
    uses discontinuous piecewise-linear basis functions, while $K_1$ uses
    continuous piecewise-linear basis functions on the same mesh.  
    In 2D, there are approximately twice as
    many elements in a regular triangular mesh as there are nodes, and three pressure
    DoFs per element for the \ponedisc{} discretization versus one per node
    for \pone{}.  Thus, there
    are approximately six %
    times as many DoFs in the \ponedisc{} pressure field than in the \pone{} pressure field.
Since $\pone{}\subset\ponedisc{}$, $P^p_0$ is
    naturally given by finite-element interpolation, determining the
    coefficients in the \ponedisc{} basis that correspond to a given function in
    \pone{}.  In the case when discontinuous pressure degrees of freedom are
    co-located
at mesh nodes (with one DoF for each adjacent element), this mapping is
simply the duplication of a continuous pressure value at the node to each
discontinuous DoF associated with the node, and its discrete structure is easy
to infer from the mesh and associated operators.  Since
$\ponedisc{}\not\subset\pone{}$, however, direct
injection or averaging %
when restricting from
\ponedisc{} to \pone{} space is not as trivial.  Here, we define
$R^p_0$ applied to a function,
$p_{dg}\in\ponedisc{}$, as yielding its corresponding continuous analog, $p_{cg}$,
by finite-element projection.
That is, we compute $p_{cg}\in\pone{}$
    that satisfies $\langle p_{cg},q_{cg}\rangle = \langle p_{dg},q_{cg}
\rangle$ for all $q_{cg}$ in the \pone{} space. As $p_{cg}$ and $q_{cg}$
can be expressed in terms of basis functions, this leads to solving a matrix
equation involving a standard \pone{} mass matrix, where the right-hand side is
defined in terms of basis functions of the two spaces.
In either case, $P_0$ and $R_0$ are assembled when %
$K_0$ and $K_1$ are assembled, as the finite-element assembly machinery for the Stokes 
operators can be repurposed for the grid-transfers.
If required, FGMRES preconditioned with a standard scalar AMG solver is used for the 
mass matrix solve, requiring a reduction  
of the $\ell_2$ norm of the relative
residual to $10^{-12}$. This is typically achieved without an increase in
iteration counts as the mesh is refined. %
While simpler
approximations of the mass-matrix inverse are certainly possible; preliminary experiments showed that these result in increased iteration counts for the overall solver.
Similar uses of AMG for mass-matrix solves have appeared
in ~\cite{geenen2009scalable}, which reports 
near-constant iteration counts as the mesh is refined, as well as in~\cite{adler2021monolithic} %
and~\cite{ur2011iterative}.

\begin{remark}[\ptwopone{} to \isoptwopone{} Grid-Transfer Operators]
  A natural question is whether %
  finite-element projection
    should be used to transfer residuals and corrections between the
    non-nested velocity spaces $\pmb{ \mathbb{P}}_2$ and $\pmb{\mathbb{P}}_1
    \text{iso}\kern1pt\pmb{ \mathbb{P}}_2$.
    To examine this in a simpler setting, we apply 
    AMG to the 
    \isoqtwo{} discretization as a preconditioner for the 
    \qtwo{} discretization of a scalar Poisson problem on the unit square.  
    Defining $R_0$ via finite-element projection 
    does not
    significantly improve performance, leading to a maximum gain of only
    one or two iterations in any given solve.  Given the 
    increased cost and complexity of this projection,
    we do not adopt %
    this for grid transfers
    between $\pmb{ \mathbb{P}}_2$ and $\pmb{\mathbb{P}}_1
    \text{iso}\kern1pt\pmb{ \mathbb{P}}_2$. However, %
    projections may likely show more benefit in other settings, e.g., 
    for even higher-order discretizations.
\end{remark}

\subsection{Smoothed Aggregation AMG for the Stokes Operator}

The standard two-level error-propagation operator for a
multigrid cycle is given by
\begin{equation}
  G_{\ell} = (I - \omega_{\ell} M_{\ell}^{-1} K_{\ell})^{\nu_2} (I - P_{\ell}
    K_{\ell+1}^{-1} R_{\ell} K_{\ell}) ( I -\omega_{\ell} M_{\ell}^{-1}
    K_{\ell})^{\nu_1} \text{ for } \ell>0, \label{eq:TG-operator}
\end{equation}
where $\ell$ denotes the level within the multigrid hierarchy,
$K_{\ell}$ is the level-$\ell$ system matrix,
$M_{\ell}$ represents the relaxation operator
applied to $K_\ell$ with damping parameter $\omega_\ell$,
$P_\ell$ is the interpolation between levels $\ell+1$ and $\ell$,
and $R_\ell$ is the restriction from level $\ell$ to level $\ell+1$.
$G_1$, for example, defines a two-level solver
within~\eqref{eq:TG-operator-full} to approximately solve the $K_1$ system.
For algebraic multigrid methods, $K_{\ell+1}$ is almost always
taken as $K_{\ell+1}=P_\ell^T K_\ell P_\ell$, and a standard multigrid
algorithm replaces the $K_{\ell+1}$ solve with a recursive application of the
two-grid solver. The multilevel cycle is fully specified once the $P_{\ell},
R_{\ell}$, $K_{\ell+1}$, and $ M_{\ell}$ are defined for all hierarchy levels. In this
paper, we take $R_\ell = P_\ell^T$ and $K_{\ell+1} = P_\ell^T K_\ell P_\ell$
within the AMG hierarchy (for $\ell \ge 1$)
and, so, only algorithms for $P_{\ell}$ and $ M_{\ell}$ need to be determined.

While popular AMG preconditioners can be directly applied to our \isoptwopone{}
system, we find that these generally do not lead to fast convergence. %
Of course, most popular AMG algorithms were originally
designed for scalar PDEs coming from elliptic equations that often give rise to
positive-definite matrices. %
While many AMG
techniques are effective on elliptic PDE systems,  %
they are not directly applicable to saddle-point systems.
For example, standard smoothed aggregation (SA) AMG %
assumes that damped Jacobi
relaxation converges and smooths errors when applied to the system
matrix, neither of which is true for saddle-point systems such as $K_1$, for
which Jacobi is not even well-defined due to the zero diagonal block in the
matrix.

To apply AMG to our system, we define block-structured grid transfer
operators, as were used in~\eqref{eq:interp}, %
and then apply
AMG algorithms to define the individual blocks of the interpolation matrix.
This structure ensures that fine-level
pressures (velocities) interpolate only from coarse-level pressures
(velocities), and not from coarse-level velocities (pressures). For the velocity
block, since the vector Laplacian is an elliptic PDE discretized using
vector $\pmb{\mathbb{P}}_1$ elements, most AMG algorithms can be reliably
applied to $A$ to produce a suitable $P^{\vec{u}}_\ell$ for the velocity
interpolation operator. AMG cannot, however, be directly applied to the zero
diagonal block associated with pressure, so some auxiliary operator is needed to
determine the pressure interpolation. One possible choice for this
auxiliary matrix is the associated pressure Laplacian, $A_p$, discretized using
the \pone{} pressure basis functions.  %
Using $A_p$ has the disadvantage that it requires an
additional matrix for the AMG setup process, albeit one that is easily generated
using any finite-element package.  An alternative and purely algebraic auxiliary
operator is $BB^T$ (or some approximation to $BM_u^{-1}B^T$, where $M_u$ is the velocity mass matrix), perhaps in conjunction with some type of stencil truncation~\cite{prokopenko2017algebraic}. $B B^T$ represents a discrete
Laplace operator, although it generally has a somewhat wider stencil than
$A_p$, which can prove challenging for AMG coarsening algorithms. In the
remainder of the paper, we use $A_p$, assembled %
on the fine grid, to avoid the introduction of additional tuning
parameters associated with truncating $BB^T$ (and in approximating the inverse of $M_u$ to form an approximation to $BM_u^{-1}B^T$).

\Cref{alg:stokes_mg_setup} details our multigrid hierarchy setup for the
2D case, assuming either direct access to each block of the Stokes systems
in~\cref{eq:saddle}, or being able to infer $A$ and $B$ from $K_1$, based on the
knowledge of the total number of discrete degrees of freedom in each field.
While a systems version of SA-AMG could be applied directly to the
vector Laplacian, we instead apply a scalar version of SA-AMG to each component
of the vector Laplacian.  As expected, the two approaches
yield similar aggregates and interpolation operators,
assuming the systems SA-AMG algorithm is provided the two-dimensional
near-kernel space,
\begin{equation}\label{eq:near_null}
    {\cal N}_v
    =
    \begin{bmatrix}
        \vec{1}_{n_{x}}&           \\
                    &  \vec{1}_{n_{y}}\\
    \end{bmatrix}.
\end{equation}
where $\vec{1}_k$ refers to a length $k$ vector with all entries equal to one.
We note that the systems approach may be advantageous when considering problems
with complex boundary conditions, or generalized operators where $A$ is not a
block-diagonal matrix. %
\begin{algorithm}
	    \begin{algorithmic}[1]
        \State \textbf{Input:}  $K_1=\begin{bmatrix} A & B^T\\ B &
        \end{bmatrix}$,  Stokes system where $B \in \mathbb{R}^{n_{p}\times
            (n_{x}+n_{y})}$
	        \State  \hspace{3.1em} $A_p$, pressure stiffness matrix or a truncated version of $BB^T$
	        \State \textbf{Output:}  $K_1, $\ldots$, K_{\ell_{max}}$, Grid hierarchy
	        \State  \hspace{4.1em} $P_1, $\ldots$, P_{\ell_{max}-1}$, Interpolation Operators
	        \State
	        \State $A_x$, $A_y$ = split(A)            \Comment{Split vector Laplacian component-wise}
          \State $\{P_1^{x},\dots,P_{\ell}^{x}, \dots
            P_{\ell^{x}_{\text{max}}-1}^{x}\} = \text{sa\_amg}(A_x,
            \vec{1}_{n_{x}})$\label{line:sax} \Comment{Build interpolation for velocity in $x$}
          \State $\{P_1^{y},\dots,P_{\ell}^{y}, \dots
            P_{\ell^{y}_{\text{max}}-1}^{y}\} = \text{sa\_amg}(A_y,
            \vec{1}_{n_{y}})$\label{line:say} \Comment{Build interpolation for velocity in $y$}
          \State $\{P_1^{p},\dots,P_{\ell}^{p}, \dots
            P_{\ell^{p}_{\text{max}}-1}^{p}\} = \text{sa\_amg}(A_p,
            \vec{1}_{n_{z}})$\label{line:sap}\Comment{Build interpolation for pressure}
	        \State $\ell_{max} = \text{min}\left(\ell^{x}_{\text{max}}, \ell^{y}_{\text{max}}, \ell^{p}_{\text{max}}\right)$  \Comment{Calculate levels in grid hierarchy}
	        \State
	        \For{$\ell = 1,\ldots,\ell_{max}-1$}
	            \State $P_\ell = \begin{bmatrix} P^{x}_\ell & & \\ & P^{y}_\ell & \\ & & P^p_\ell \end{bmatrix}$ 	\Comment{Assemble monolithic interpolation operator}
	            \State $K_{\ell+1} = P_{\ell}^T K_{\ell} P_{\ell}$ \Comment{Compute coarse-grid operator}
	        \EndFor%
	    \end{algorithmic}
	    \caption{SA-AMG Hierarchy Setup for 2D Stokes System}\label{alg:stokes_mg_setup}
\end{algorithm}

Since we independently coarsen the velocity and pressure blocks, SA-AMG may
determine different numbers of levels in each hierarchy.  Thus, we
form a block hierarchy using the fewest levels across these
hierarchies.  While this could lead to impractically large coarsest-grid systems
if the coarsening of each field is dramatically different, we have yet to encounter
any such problems. After the interpolation
operators for each field are formed, they are combined
as a block-diagonal matrix into a single monolithic interpolation operator
$P_{\ell}$, and used in a monolithic calculation of the coarse-grid operator.
We emphasize that we use smoothed aggregation completely independently on each
field (lines \ref{line:sax},\ref{line:say} and \ref{line:sap} in~\cref{alg:stokes_mg_setup}), without
any explicit coordination or coupling (in contrast to most traditional monolithic MG
methods). As a result, the $B$ block in $K$ is projected to its coarse-grid
analog using both the velocity and pressure grid transfers.
A natural concern is whether this independent coarsening might lead to somewhat
strange stencils in $B$, where the coarse-grid discretization might not
necessarily even be stable, even if the fine-level $B$ is stable.
While our approach does not guarantee stability, we have observed that we
generally get stable coarse-grid operators so long as we maintain a similar
ratio between the velocity and pressure DoFs throughout the hierarchy, which is
crucial to the convergence of the overall AMG method.
This can be achieved when we (approximately) match the coarsening rates for each
field on each level. These coarsening rates are effectively determined by the
strength-of-connection (SoC) measures used in aggregation, which identify the
edges of the matrix graph that can be ignored during the coarsening
process.  Using traditional SoC measures, such as the symmetric strength of
connection commonly used in smoothed
aggregation~\cite{PVanek_JMandel_MBrezina_1996a}, we found that achieving
compatible coarsening rates was difficult, requiring significant tuning of
drop-tolerance parameters for the velocity and pressure fields at each level
of the hierarchy. Instead, we employ the so-called \textit{evolution} SoC
from~\cite{olson2010new}, which provides more uniform coarsening by integrating
local representations of algebraically smooth error into the SoC measures.  We
observe (experimentally) that this leads to consistent coarsening rates across
the velocity and pressure fields, with no need for parameter
tuning~\cite{proceedings21}.

\begin{remark}[Applying~\cref{alg:stokes_mg_setup} directly to $K_0$]\label{rmrk:ho_solvers}
\cref{alg:stokes_mg_setup} can be applied directly to
the \ptwopone{} discretization, %
resulting in a
preconditioner that we denote by \Mh{}. While the effectiveness of \Mh{} is
expected to decline for discretizations with higher-order bases due to
mismatch in pressure and velocity field coarsening rates, we find that it offers
stable performance for the \ptwopone{} discretization. %
However, direct application of~\cref{alg:stokes_mg_setup} to the \ptwoponedisc{}
discretization leads to poorly convergent or
non-convergent preconditioners, which we attribute to limitations
in the algebraic coarsening of the discontinuous pressure field.
\end{remark}

\subsection{Relaxation}\label{sec:mg_relax}

The final piece of the AMG method is the %
relaxation method with (weighted) error propagation operator $I - \omega_\ell
M_\ell^{-1}K_\ell$. 

\subsubsection{Additive Vanka}\label{sec:vanka_relax}

Vanka relaxation~\cite{SPVanka_1986a} is a standard relaxation 
scheme for saddle-point problems. %
It incorporates coupling between PDE components using local ``patch'' solves over small sets of DoFs (similar to small sub-domains within an additive or 
multiplicative Schwarz method).

A Vanka patch is typically constructed to contain a small number of
pressure DoFs (often just one) and all velocity DoFs that are adjacent (i.e.,
have nonzero connections) to this ``seed'' pressure DoF in the system matrix,
$K$.
Within geometric multigrid, such patches can be naturally constructed
using topological arguments~\cite{farrell2019c}; here, however, we construct the
patches algebraically on each level of the hierarchy.
Although the general
mechanism for patch assembly is the same for all discretizations, we define
each patch with a single seed pressure DoF for both the \isoptwopone{}
multigrid hierarchy and the fine-grid \ptwopone{} discretization. In these cases,
each row of $B$ corresponds to one patch, and patch $i$ incorporates
DoF $i$ of the pressure and all velocity DoFs linked to nonzero entries in row
$i$ of $B$. %

\begin{figure}[!ht]
  \centering
  \begin{subfigure}{0.31\textwidth}
    \includegraphics[width=\textwidth]{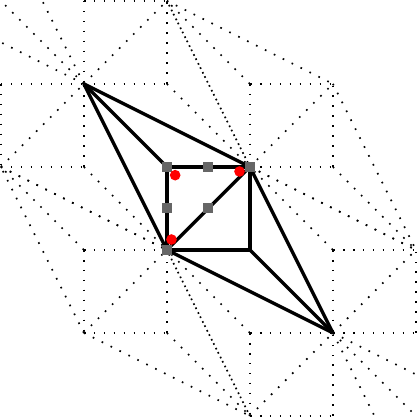}
    \caption{\ptwoponedisc{}}
  \end{subfigure}
  \hfill
  \begin{subfigure}{0.31\textwidth}
    \includegraphics[width=\textwidth]{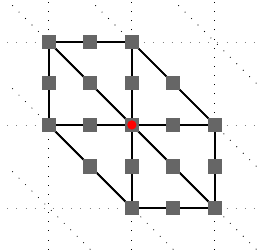}
    \caption{\ptwopone{}}
  \end{subfigure}
  \hfill
  \begin{subfigure}{0.31\textwidth}
    \includegraphics[width=\textwidth]{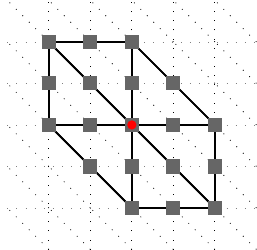}
    \caption{\isoptwopone{}}
  \end{subfigure}
  \captionsetup{skip=1pt}
  \caption{(Left) one Vanka patch for the \ptwoponedisc{} discretization;
    (middle) one Vanka patch for the \ptwopone{} discretization; and (right) one
    Vanka patch for the \isoptwopone{} discretization.
    The pressure DoF(s)
    (red circles, \protect\tikz[baseline=-0.5ex]\protect\draw[red, fill=red]
    (0, 0) circle (2pt);)
    and the velocity DoFs directly algebraically connected to that pressure DoF
    (gray squares, \protect\tikz[baseline=-0.5ex]\protect\fill[black!50,
    fill=black!50] (-2pt, -2pt) rectangle (2pt, 2pt);)
    are marked. Displayed patches are on 2D triangular meshes.
}\label{fig:Vanka_patches}
\end{figure}

For the \ptwoponedisc{} discretization, we employ three seed pressure DoFs for
each patch, corresponding to the three pressure DoFs within each element.
The three pressure DoFs on each element can be identified by three rows of $B$
with identical sparsity patterns, making it possible to determine the seeds
algebraically.
However, since the information is readily available and requires minimal
pre-processing, we opt to use a straightforward element-to-DoF map to identify
the three pressure DoFs on a single element.
\Cref{fig:Vanka_patches} depicts typical patches for the \ptwoponedisc{},
\ptwopone{} and \isoptwopone{} discretizations of the fine-level operators.
While these patches are directly linked to element/mesh entities on the finest
mesh, a clear element/mesh interpretation is absent on algebraically
coarsened levels.

On level $\ell$, the $i^{th}$ patch matrix $K_\ell^{(i)}$ is then given by a
square sub-matrix, which is composed of patch sub-blocks $B^{(i)}$
and $A^{(i)}$, written as
\begin{equation}%
    K_\ell^{(i)}
    =
    \begin{bmatrix}
        A^{(i)}  & (B^{(i)})^T \\
        B^{(i)}  &  0
    \end{bmatrix}.
\end{equation}
Formally, $K_\ell^{(i)} =V_\ell^{(i)} K_\ell (V_\ell^{(i)})^T$ where
$V_\ell^{(i)}$ is an injection operator defined for each patch, $i$, on level
$\ell$ as an $m^{(i)}_\ell \times m_\ell$ binary matrix that extracts the
subset of global DoFs on level $\ell$ coinciding with the $i^{\text{th}}$ patch, where $m_\ell$ is
the total number of rows/columns of $K_{\ell}$, while $m^{(i)}_\ell$ is the number of DoFs
contained in $K_\ell^{(i)}$, which is one more than the number of velocity DoFs
in the $i^{th}$ patch for the \isoptwopone{} and \ptwopone{} discretizations,
and three more for the \ptwoponedisc{} discretization.

In this work, we use the additive form of Vanka relaxation~\cite{schoberl2003schwarz,
farrell2021local}. A single Vanka sweep for \isoptwopone{} and \ptwopone{} discretizations 
is given by
\begin{equation*}\label{eq:vanka}
  M^{-1}_{\ell} = \sum_{i=1}^{n^{(p)}_{\ell}} (V_\ell^{(i)})^{T} W_\ell^{(i)} \left (V_\ell^{(i)} K_\ell (V_\ell^{(i)})^T \right )^{-1}  V_\ell^{(i)},
\end{equation*}
where the diagonal weighting matrix, $W_\ell^{(i)}$, is a partition of unity, with each
diagonal entry equal to the reciprocal of the number of patches that contain
the associated degree of freedom.  For the \ptwoponedisc{} discretization, the
sum extends over fewer patches, for $1\leq i \leq n^{(p)}_{\ell}/3$.

Applying inverses of patch matrices often constitutes the most
computationally demanding part of the multigrid method. %
The direct inverse of the patch can be applied using highly
optimized dense-linear algebra kernels. However, this requires
significant storage, as it results in a fully dense inverse matrix over each patch.
A sparse LU factorization reduces the memory footprint but
results in slower solution times due to sparse forward and backward
substitutions.
Our performance studies, not presented here, show that sparse LU does not
outperform the direct-inverse method in efficiency. However, the dense inverses'
high storage may be limiting.

Therefore, we introduce a novel approach that balances the effectiveness of
dense-linear algebra with reduced storage for sparse LU
factorizations. Instead of computing the inverse of each $K_\ell^{(i)}$, we
leverage a block LU factorization of each Vanka patch,
\begin{equation}\label{eq:vanka_fact1}
     K_\ell^{(i)} =
       \begin{bmatrix}
          A  & B^T \\
          B   &  0
      \end{bmatrix}
      =
      \begin{bmatrix}
          A   &  \\
          B   &  S
      \end{bmatrix}
      \begin{bmatrix}
          I  & U  \\
             & I
      \end{bmatrix},
\end{equation}
with sub- and superscripts on the sub-blocks omitted for clarity. Here, 
$U=A^{-1} B$, and $S=-BA^{-1}B^T$ denotes the Schur complement.  We
emphasize that $B^T$ and $U$ have small column dimension (one in the
typical case, but three for the \ptwoponedisc\ discretization), and that $S$ is
a small matrix.  The action of the inverse for $K_\ell^{(i)}$
is 
\begin{equation} \label{eq:vanka_fact2}
     K_\ell^{(i)}  x_\ell^{(i)}=
     b_\ell^{(i)} \; \Rightarrow \;
     x_\ell^{(i)} =
      \begin{bmatrix}
          I  & \hat{U} \\
             & I
      \end{bmatrix}
      \begin{bmatrix}
          A^{-1}  &  \\
          \hat{B}  &  S^{-1}
      \end{bmatrix}
      b_\ell^{(i)}.
\end{equation}
where $\hat{B}=-S^{-1}B A^{-1}$, $\hat{U} = -U$, and $\hat{B}$, $\hat{U}$, and $S^{-1}$ are 
precomputed in the setup phase.
We further leverage the block-diagonal nature of the vector Laplacian
operator to minimize the storage associated with the inverse of sub-block
$A$ up to a factor of 2 and 3 in 2D and 3D, respectively. Furthermore, this method
reduces the number of floating-point operations necessary for calculating the
action of the patch's inverse, thereby enhancing the computational efficiency of
the setup and solution phases.

Vanka relaxation
exhibits sensitivity to its damping parameter, $\omega_\ell$. 
Fourier analysis is often used to determine $\omega_\ell$ when Vanka is
employed within a GMG
scheme~\cite{he2019local,farrell2021local,SPMacLachlan_CWOosterlee_2011a,LFAoptAlg2020};
however, this method is not feasible for the irregular coarse meshes generated in an
AMG algorithm.  Several recent papers have used alternative
methods to tune relaxation parameters.  In~\cite{adler2021monolithic}, for
example, instead of a stationary weighted relaxation scheme, a Chebyshev
iteration is used, with hand-tuned Chebyshev intervals chosen to
optimize the resulting solver performance.  While this approach is, in
principle, also applicable here, it could require significant tuning for the
various levels of the multigrid hierarchy.  Instead, we opt for a more practical
approach where relaxation on each level is accelerated using FGMRES\footnote{We use FGMRES instead of classical GMRES to avoid the extra application of a (possibly expensive) preconditioner needed by classical GMRES to recover the approximate solution.}.
FGMRES does not necessarily provide optimal relaxation weights as it constructs
a residual-minimizing polynomial, which is not the same as optimally damping
high-frequency error.
While  it does not provide general convergence bounds as can sometimes be achieved using
Chebyshev iterations, it does eliminate the need for additional tuning
parameters and generally provides satisfactory damping weights.
One minor drawback is that FGMRES requires additional inner products,
though this is a minimal computational cost over the Vanka relaxation patch solves.

In the remainder of the paper, we fix the relaxation to include two inner FGMRES
iterations with one Vanka relaxation sweep per iteration, akin to
 setting $\nu_1=\nu_2=2$ in~\eqref{eq:TG-operator-full} or~\eqref{eq:TG-operator}.
Results below show that using this Krylov wrapping approach on all levels in the hierarchy provides
stable and robust convergence for systems discretized with \ptwopone{} elements.
However, we observed stagnant convergence when applied to systems discretized
with \ptwoponedisc{} elements. To address this issue, we instead use a
stationary Vanka iteration on the finest grid with the \ptwoponedisc{} elements
in~\eqref{eq:TG-operator-full} and Krylov acceleration on all other levels
in~\eqref{eq:TG-operator}.  We perform a parameter search to identify the
optimal $\omega_0$ value, denoted as $\omega_0^{opt}$, for the $K_0$ systems
discretized with the \ptwoponedisc{} element pair, discussed below.

\subsubsection{Distributive least squares commutator (DLSC)}\label{sec:dlsc}

We consider a distributive least squares commutator (DLSC) method that builds
on the least-squares commutator distributive Gauss-Seidel (LSC-DGS) approach
outlined in~\cite{wang2013multigrid}. LSC-DGS has shown to be effective for
various Stokes problems, including its adaptation to linearized Navier-Stokes
equations~\cite{elman2006block} and as a defect-correction preconditioner in
low-order Stokes problems~\cite{wang2013multigrid}. The LSC-DGS method is
derived from Schur complement approximations using a least-squares-based
approximate commutator. The notable advantage is the elimination of
user-defined parameters, unlike in the classical DGS method. For an in-depth study of the
LSC-DGS method in the context of the Marker and Cell (MAC) discretization, readers are referred to~\cite{wang2013multigrid}.

DLSC methods rely on the commutative assumption that, under certain conditions, $A M_u^{-1} B^T \approx B^T M_p^{-1} A_p$, where $A_p$ represents a discrete Laplacian on the pressure, and $M_u$ and $M_p$ are velocity and pressure mass-matrices, respectively. This assumption leads to an approximate factorization, which can be expressed as:
\begin{subequations}
\begin{align*}
&\begin{bmatrix}
A & B^T\\
B & 0
\end{bmatrix}
\begin{bmatrix}
M^{-1}_u & M^{-1}_u B^T\\
0 & -M_p^{-1}A_p
\end{bmatrix}
\approx
\begin{bmatrix}
AM^{-1}_u & 0\\
BM^{-1}_u & B M^{-1}_u B^T
\end{bmatrix}, \ \\
\text{or}\quad
&\begin{bmatrix}
A & B^T\\
B & 0
\end{bmatrix}^{-1}
\approx
\begin{bmatrix}
M^{-1}_u & M^{-1}_u B^T\\
0 & -M_p^{-1}A_p
\end{bmatrix}
\begin{bmatrix}
AM^{-1}_u & 0\\
BM^{-1}_u & B M^{-1}_u B^T
\end{bmatrix}^{-1}.
\end{align*}
\end{subequations}
Further simplifications are achievable by setting $$A_p = M_p\hat{A}_p = M_p (BM^{-1}_u B^T)^{-1}(BM^{-1}_u AM^{-1}_u B^T)$$ and factorizing the lower-triangular matrix:
\begin{subequations}
\begin{align*}
\begin{bmatrix}
M^{-1}_u & M^{-1}_u B^T\\
0 & -M_p^{-1}(M_p \hat{A}_p)
\end{bmatrix}
\begin{bmatrix}
M_u & \\
 & I
\end{bmatrix}
&\begin{bmatrix}
A & 0\\
B & B M^{-1}_u B^T
\end{bmatrix}^{-1}\\
\begin{bmatrix}
I & M^{-1}_u B^T\\
0 & -\hat{A}_p
\end{bmatrix}
&\begin{bmatrix}
A & 0\\
B & B M^{-1}_u B^T
\end{bmatrix}^{-1}.
\end{align*}
\end{subequations}

The DLSC iteration for the linear system
\begin{equation*}
\begin{bmatrix}
A & B^T\\
B & 0
\end{bmatrix}
  \begin{bmatrix}
\vec{u}\\ 
p
\end{bmatrix}
=
\begin{bmatrix}
\vec{f}\\ 
g
\end{bmatrix}
\end{equation*}
is then given by
\begin{equation}\label{eq:fulldlsc}
\begin{bmatrix}
\vec{u}\\ 
p
\end{bmatrix}
\leftarrow
\begin{bmatrix}
\vec{u}\\ 
p
\end{bmatrix}
+
\begin{bmatrix}
I & M^{-1}_u B^T\\
0 & \hat{A}_p
\end{bmatrix}
\begin{bmatrix}
A & 0\\
B & B M^{-1}_u B^T
\end{bmatrix}^{-1}
\begin{bmatrix}
\vec{f} - A \vec{u} - B^T p\\
g - B \vec{u}
\end{bmatrix}.
\end{equation}
Notably, operator  $M_p$ does not appear in the final form of the algorithm as seen in \cref{eq:fulldlsc} and \cref{alg:DLSC}. 
In the context of this work, unless stated otherwise, $M_u$ is replaced by the diagonal of the mass matrix. This mass matrix scaling is consistently applied across all \ptwopone{} grids. For the coarse-level velocity mass matrix, we utilize the Galerkin product, constructed using the grid-transfer operators from~\cref{alg:stokes_mg_setup}. On \isoptwopone{} grids, the velocity mass matrix is approximated by an identity operator, as using a diagonal mass-matrix approximation does not meaningfully impact the convergence.

\Cref{alg:DLSC} describes a single sweep of this relaxation, where the first two solve steps
(lines~\ref{line:momentum} and \ref{line:continuity}) result from the block lower
triangular forward solve and where the last solve and correction
(line~\ref{line:pressure_correction} and line~\ref{line:velocity_correction}) results from the block upper triangular
matrix multiplication in~\cref{eq:fulldlsc}. A central aspect in the LSC-DGS
algorithm of~\cite{wang2013multigrid} is the use of Gauss-Seidel
to approximate the solutions to the $A$ and $BM_u^{-1}B^T$ matrix systems.  In
contrast, we make use of Chebyshev polynomial relaxation (and, later, a standalone AMG V-cycle), 
which leads to improved parallel efficiency over Gauss-Seidel.
Our empirical observations for DLSC, which are largely in agreement with the
findings in~\cite{wang2013multigrid}, reveal that the pressure-correction step
(line~\ref{line:pressure_correction}) typically requires more iterations compared to the
momentum relaxation step (line~\ref{line:momentum}) and the continuity
relaxation step (line~\ref{line:continuity}). This requirement becomes
increasingly pronounced in the context of anisotropic unstructured meshes and
three-dimensional problems. 
This loosely translates to the fine-grid velocity being closer to being divergence-free when the pressure-correction step is solved more accurately.
\begin{algorithm}
  \begin{algorithmic}[1]
\State \textbf{Input:} $\vec{u}$, $p$, initial guess
\State \hspace{3.1em}  $\vec{f}$, $g$, components of the right-hand side
\State \textbf{Output:} $\vec{u}^{\text{new}}$, $p^{\text{new}}$, updated variables \\
\State
$\vec{v}_*\leftarrow \text{Approx. solve}\, A \vec{v}_* = \vec{f} - A \vec{u} - B^T p$\Comment{solve the momentum equation}\label{line:momentum}
\State $q\kern5pt\leftarrow \text{Approx. solve}\, B M^{-1}_u B^T q = g -  B(\vec{u}+\vec{v}_*)$ \Comment{solve the continuity equation}\label{line:continuity}
\State 
\State $\delta p \leftarrow \text{Approx. solve}\, B M^{-1}_u B^T \delta p =
      -(BM^{-1}_u AM^{-1}_u B^T)q$\Comment{pressure correction}\label{line:pressure_correction}
\State $\delta\vec{u} =  \vec{v}_* + M^{-1}_u B^T q$\Comment{velocity correction}\label{line:velocity_correction}
\State 
\State $\vec{u}^{\text{new}} = \vec{u} +\delta\vec{u}$ \Comment{velocity update}
\State $p^{\text{new}} = p + \delta p$ \Comment{pressure update}
  \end{algorithmic}
  \caption{DLSC} \label{alg:DLSC}
\end{algorithm}

Numerical results in~\cref{sec:results} show that this is an effective relaxation scheme for the \ptwopone{} and \isoptwopone{} discretizations, but we did not find it to be effective for the \ptwoponedisc{} discretization.  
We note that~\cite{wang2013multigrid} proposes the use of element-wise block
Gauss-Seidel relaxation for $BM^{-1}_u B^T$ in a similar setting; block variants such as
block-Jacobi-Chebychev, either standalone or within an AMG V-cycle, may also be
effective, but we leave this open question to future work. In all SV results
in~\cref{sec:results}, we only use cell-wise Vanka relaxation (and not DLSC)
for the \ptwoponedisc{} discretization.

\section{Numerical Results}\label{sec:results}

\begin{figure}[!ht]
    \centering
    \includegraphics[width=\textwidth]{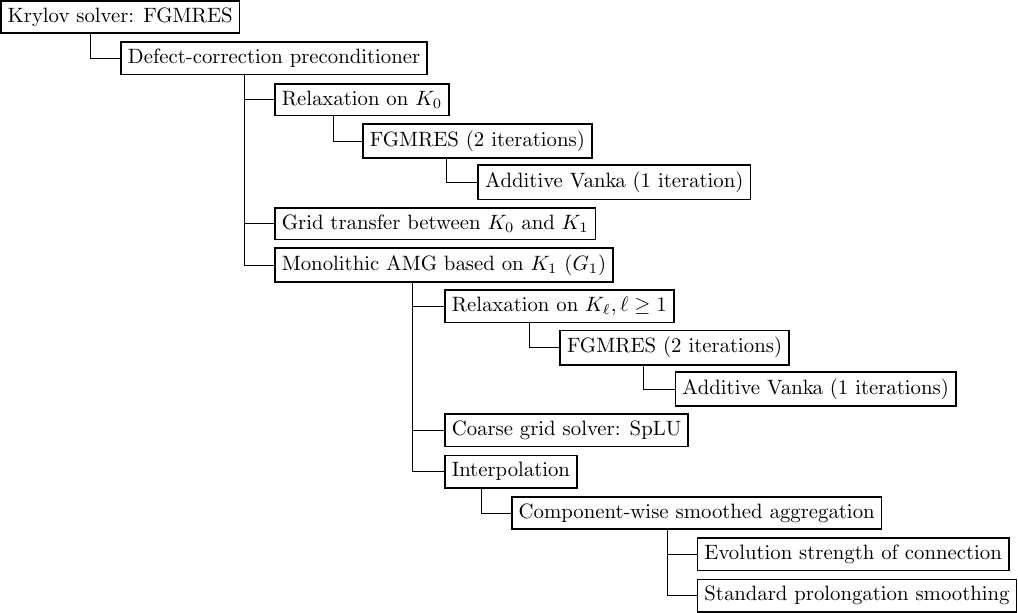}
    \caption{Solver diagram for defect-correction preconditioners. If relaxation
    is performed on all levels, the diagram corresponds to the \cMphHLO{}
    solver.}\label{fig:solver_diagram}
\end{figure}

We compare the \Mh{} preconditioner with the proposed defect-correction 
AMG preconditioner using two different relaxation types: additive Vanka and DLSC.
The \Mh{} preconditioner is
constructed by applying~\cref{alg:stokes_mg_setup} directly to the high-order
discretization ($K_0$), while defect correction %
applies~\cref{alg:stokes_mg_setup} to an auxiliary operator, $K_1$.
For the defect-correction (DC) preconditioners, we consider three different
relaxation schedules: \cMphHLO{}, \cMphHO{}, and \cMphLO{}.
The \cMphHLO{} preconditioner performs pre- and post-relaxation on all levels of the hierarchy,
except for the coarsest level ($n_{lvls}$) where an exact solution is used.
\cMphHO{} skips relaxation on the $K_1$ system, but
performs relaxation on the $K_0$ system and remaining levels of the
\isoptwopone{} hierarchy, except for a direct solve on the coarsest level.
\cMphLO{} skips relaxation on the $K_0$ system and performs pre- and
post-relaxation on remaining levels, except for a direct solve on the coarsest level.
A solver diagram for the defect-correction approach using Vanka relaxation is given
in~\cref{fig:solver_diagram}. \Mh{} is similar (incorporating our
AMG enhancements) except that the Krylov solver is directly applied to monolithic
AMG based on  $K_0$.
We will also explore additional preconditioning variants (described later) based on the defect-correction
operator. The notation for all relevant parameters and
solvers is summarized in~\cref{table:param_ref}.

\begin{table}[!ht]
\centering
\begin{tabular}{cl}
\toprule
   Symbol &Description\\
\midrule
$K_0$     & saddle-point system for the higher-order discretizations (\ptwopone{}, \ptwoponedisc{}) \\
$K_1$     & saddle-point system for the lower-order discretization (\isoptwopone{}) \\
$G_1$     & monolithic AMG Solver based on $K_1$ \\
\midrule
$\ell$    & generic level in the monolithic multigrid hierarchy \\
$n_{lvls}$& total number of levels/grids in the monolithic multigrid hierarchy \\
$\eta_u$  & momentum equation residual scaling parameter on level $0$ \\
$\eta_p$  & continuity equation residual scaling parameter on level $0$ \\
$\eta$    & $K_0$ residual scaling parameters when $\eta_u=\eta_p$ \\
    $\omega_{\ell}$& Vanka relaxation damping parameter \\
Superscript $opt$& denotes that the optimal parameter was computed via brute-force
    search \\
\midrule
\cMBTP{} &  Inexact block triangular preconditioner based on $K_0$ \\
    \cMh{}     &  Monolithic AMG algorithm applied directly on $K_0$\\
\cMphHLO{} & DC with relaxation on all levels,
    $\omega_{\ell}=1~\text{for}~\ell\in[0,n_{lvls})$ \\           %
\cMphLO{}  & DC with relaxation on all but the $K_0$ level,
    $\omega_{\ell}=1~\text{for}~\ell\in [1,n_{lvls})$ \\          %
\cMphHO{}  & DC with relaxation on all but the $K_1$ level,
    $\omega_{\ell}=1~\text{for}~\ell\in [0]\cup[2,n_{lvls})$ \\   %
\midrule
\tikz[baseline=-0.5ex]\draw[solid] (0,0) -- (1,0); Line & DC preconditioner with Vanka relaxation on~$\ell\in[0,n_{lvls})$ \\
\tikz[baseline=-0.5ex]\draw[dotted] (0,0) -- (1,0); Line& DC preconditioner with DLSC relaxation on~$\ell\in[0,n_{lvls})$\\
\tikz[baseline=-0.5ex]\draw[dash dot] (0,0) -- (1,0); Line & DC preconditioner with Vanka on $\ell=0$ and DLSC on $\ell\in[1,n_{lvls})$ \\
\bottomrule
\end{tabular}
  \caption{Notation Reference Table}\label{table:param_ref}
\end{table}

In the defect-correction preconditioners using DLSC, we perform two
iterations of DLSC without Krylov acceleration during relaxation. The optimal
degree of the Chebyshev polynomial depends on the stage of~\cref{alg:DLSC}
and the mesh type, as outlined in~\cref{table:lsd_dgs_config}. 
For 3D meshes, higher Chebyshev polynomial degrees are necessary to achieve convergence
rates similar to Vanka relaxation. 
We note that the pressure-correction step is
often the limiting factor for convergence, particularly for 3D
problems, where higher polynomial degrees are needed to achieve convergence
rates comparable to 2D cases.
\begin{table}[!ht]
\begin{center}
    \begin{tabular}{l | c c}
    \toprule
    Update Stage & 2D & 3D \\
    \midrule
      Momentum equation,   \cref{alg:DLSC}~line~\ref{line:momentum}    & 3 & 4 \\
      Continuity equation, ~\cref{alg:DLSC}~line~\ref{line:continuity}  & 3 & 4 \\
    Pressure correction, ~~\cref{alg:DLSC}~line~\ref{line:pressure_correction}   & 6 & 16 \\
    \bottomrule
    \end{tabular}
\end{center}
\caption{Degrees of Chebyshev polynomials used for polynomial relaxation at
different stages of the DLSC algorithm (\cref{alg:DLSC}). These degrees vary
depending on the dimension (2D or 3D) of the problem. Chebyshev polynomials on the interval
$[0.5^{d}\lambda, 1.1\lambda]$ are used, where $d$ represents the
problem's dimension and $\lambda$ is an estimate of the largest eigenvalue of $A$ or B$M_u^{-1}B^T$ generated by an Arnoldi iteration.  }\label{table:lsd_dgs_config} \end{table}

As a point of comparison for the monolithic AMG preconditioners, we consider a block-triangular (or inexact Uzawa) preconditioner,
denoted by \cMBTP{}. This preconditioner approximates the lower-triangular
factor in the block factorization of the Stokes
operator~\cite{elman2014finite,notay2019convergence}. It computes the velocity
and pressure updates, $[\vec{\delta u}, \delta p]$, as follows:
\begin{align*}
    \vec{\delta u} &= Q_{A}^{-1} \vec{r_u} \\
    \delta p &= Q_{B}^{-1}(B \vec{\delta u} - r_p),
\end{align*}
where $[\vec{r_u}, r_p]=b-Kx$, $Q_{A}^{-1}$ represents an approximation of the
vector Laplacian inverse, and $Q_{B}^{-1}$ denotes an exact or approximate
inverse of the pressure-mass matrix, depending on the context.
To approximate $Q_{A}^{-1}$, we employ a single V-cycle of SA-AMG\@. For
$Q_{B}^{-1}$, we use either an exact or approximate mass-matrix inverse, as
outlined in~\cite{elman2014finite}. The approach to $Q_{B}^{-1}$ is influenced
by the discretization used.
In the \ptwoponedisc{} case, $Q_{B}$ is block-diagonal, enabling us to compute
the exact inverse at low cost. However, when using the \ptwopone{}
discretization, we approximate the action of $Q_{B}^{-1}$ using a single V-cycle
of SA-AMG\@.
The SA-AMG solvers are assembled using evolution SoC, standard aggregation, and
standard prolongator smoothing for interpolation. For multigrid relaxation, we
employ two steps of Krylov-wrapped Jacobi relaxation.

The \ptwopone{}, \ptwoponedisc{}, and \isoptwopone{} systems are assembled in
Firedrake~\cite{FiredrakeUserManual,kirby2018solver}. Although
Firedrake currently does not provide the capability to construct \isoptwopone{}
matrices directly, we form them by assembling the \ponepone{}
discretization on a refined mesh followed by geometric restriction to
a coarsened pressure field. This indirect approach to forming the \isoptwopone{}
discretization can be avoided in a finite-element code base that supports
macro-elements. %
In our open-source implementation (\href{\githuburl}{\githuburltext}),
we use the PyAMG library~\cite{BeOlSc2022} to
implement AMG preconditioners that combine~\cref{alg:stokes_mg_setup} with
Vanka relaxation methods.

\begin{figure}[!ht]
  \centering
  \begin{subfigure}{0.4\textwidth}
    \centering
    \includegraphics[width=0.8\textwidth]{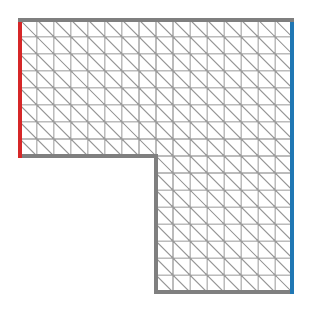}
    \caption{2D Structured}\label{fig:meshes_struct_2D}
  \end{subfigure}
  \hspace{1cm}
  \begin{subfigure}{0.4\textwidth}
    \centering
    \includegraphics[width=0.8\textwidth]{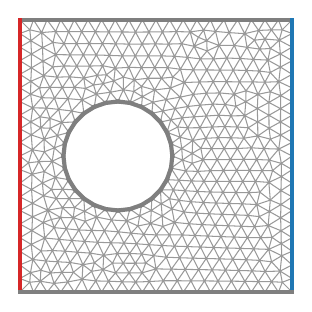}
    \caption{2D Unstructured}\label{fig:meshes_struct_3D}
  \end{subfigure}
  \\
  \begin{subfigure}{0.4\textwidth}
    \centering
    \includegraphics[width=0.7\textwidth]{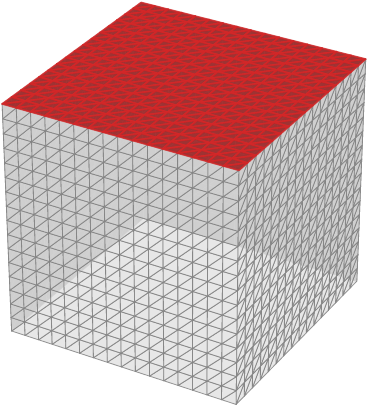}
    \caption{3D Structured}\label{fig:meshes_unstruct_2D}
  \end{subfigure}
  \hspace{1cm}
  \begin{subfigure}{0.4\textwidth}
    \centering
    \includegraphics[width=0.9\textwidth]{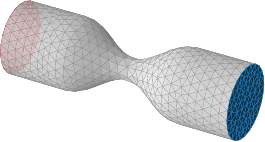}
    \vspace{1cm}
    \caption{3D Unstructured}\label{fig:meshes_unstruct_3D}
  \end{subfigure}
  \captionsetup{skip=1pt}
  \caption{Meshes used to construct the Stokes problems. The boundaries are
    marked with three colors (red, blue, and gray), corresponding to three
    types of velocity field boundary conditions (BCs). Red corresponds to non-zero
    Dirichlet boundaries. Blue corresponds to Neumann boundary conditions. Gray
    indicates zero Dirichlet BCs.}\label{fig:problem_meshes}
\end{figure}

We consider both structured and unstructured
meshes in 2D and 3D for our numerical solvers study. The structured meshes are
generated using the Firedrake
software, while the unstructured meshes are generated using
Gmsh~\cite{geuzaine2009gmsh}. Barycentric refinement of
meshes uses the code associated with~\cite{farrell2021reynolds}.  As test
problems, we consider the following, with domains pictured in
Figure~\ref{fig:problem_meshes}.
\begin{description}
\item[Structured 2D\@: Backward-Facing Step] \hfill \\
    A standard 2D
    backward-facing step domain, discretized with structured triangular
    elements~\cite{elman2014finite}. The left domain surface 
    is marked in red in Figure~\ref{fig:problem_meshes}, where a
    parabolic inflow boundary condition (BC) is imposed. The right domain surface,
    marked in blue, has a natural (Neumann) BC.  All
    other domain surfaces (marked in grey) have zero Dirichlet BCs
    on the velocity field.
\item[Structured 3D\@: Lid-Driven Cavity] \hfill \\
    A standard lid-driven cavity
    problem on a 3D domain, discretized with structured tetrahedral elements~\cite{elman2014finite}.
    The top face of the domain, marked in red, has a Dirichlet BC imposed with
        constant tangential velocity, while all other faces have zero Dirichlet
        BCs on the velocity field.
\item[Unstructured 2D\@: Flow Past a Cylinder] \hfill \\
    This test problem simulates flow past a cylinder on a 2D domain,
    discretized with unstructured triangular elements. The left 
    domain surface, marked in red, has a parabolic inflow BC imposed, while the right surface,
        marked in blue, has a natural outflow BC.  All other boundaries
        including those on the interior boundary of the cylinder, have zero Dirichlet BCs
     on the velocity field.
\item[Unstructured 3D\@: Pinched Channel] \hfill \\
    This test problem simulates flow through a pinched tube on a 3D domain,
    discretized with unstructured tetrahedral elements. The left face of the
    domain, marked in red, has a constant velocity inflow
    BC, while the right face, marked in blue, has a natural
    outflow BC. All other faces of the domain have zero
    Dirichlet BCs on the velocity field.
\end{description}

For each problem, we present iteration counts, relative solution times, 
and relative per-iteration times. 
Additionally, we include the total relative times, including the setup costs, for the two discretization types (see~\cref{fig:th_total_timings} and \cref{fig:sv_results}).
The iteration counts are
the total number of preconditioned FGMRES iterations required to reduce the
$\ell_2$ norm of the residual by a relative factor of $10^{10}$.
For all tests, we use restarted FGMRES with a maximum Krylov subspace size of $20$.
The timings are measured relative to the \cMBTP{} preconditioner.
We report the shortest of five independent time-to-solution measurements to
account for time measurement variability.

To find the best damping parameters, $\eta$ and $\omega_0$, 
we use a brute-force
search for each problem type and multigrid variant. For the \ptwopone{}
discretization, we only find $\eta$, because Krylov-wrapping
avoids the use of $\omega_{\ell}$. %
However, for the \ptwoponedisc{} discretization, Krylov-wrapping on $\ell=0$ is
ineffective. Therefore, for the \ptwoponedisc{}
\cMphHLO{} preconditioner, we calculate both $\eta$ and $\omega_0$.
We define the optimal parameter set as the value(s) that result in the fewest
number of preconditioned FGMRES iterations required to reach a relative
reduction in the $l^2$-norm of the residual by a factor
of $10^8$. All parameter searches were performed on
problems with at least one million DoFs and 4 (3) levels in the $G_{1}$ multigrid
hierarchy in 2D (3D).
Our results reveal that even a limited search yields near-optimal parameters,
indicating low sensitivity to small parameter variations. 
In cases where multiple sets of parameters result in the same number of iterations,
we select the optimal set of parameters to be the one that has the smallest
measured convergence factor,  determined by finding the
geometric mean of the residual reduction factors during the last 10 FGMRES
iterations. In the cases where the optimal parameter range includes value $1.0$,
we select it as the optimal parameter value for simplicity.

\subsection{Preconditioning \texorpdfstring{\ptwopone}{P2/P1} Systems}\label{sec:p2p1_results}

\begin{figure}[!ht]
    \centering
    \begin{subfigure}{\textwidth}
        \centering
        \includegraphics[width=\textwidth]{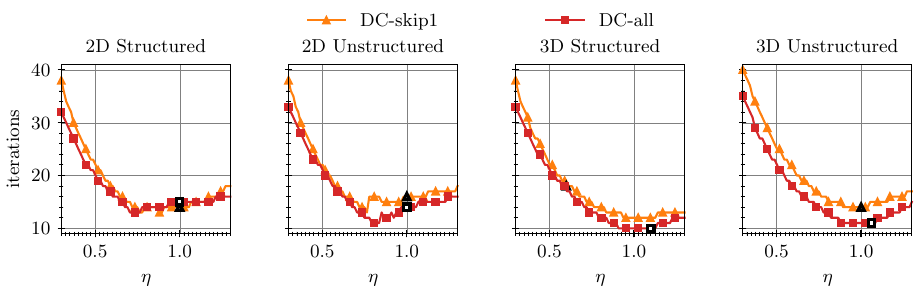}
        \captionsetup{skip=1pt}
        \caption{Sensitivity analysis for AMG hierarchies using Vanka relaxation.}\label{fig:th_sensitivity_vanka}
    \end{subfigure}
    \begin{subfigure}{\textwidth}
        \centering
        \includegraphics[width=\textwidth]{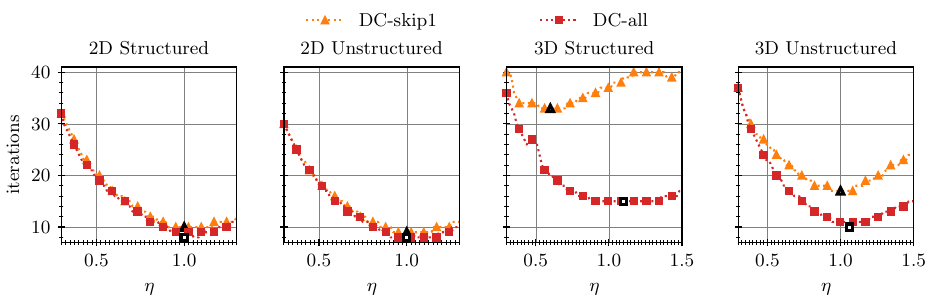}
        \captionsetup{skip=1pt}
        \caption{Sensitivity analysis for AMG hierarchies using DLSC relaxation.}\label{fig:th_sensitivity_lsc_dc}
    \end{subfigure}
    \captionsetup{skip=0.1pt}
    \caption{
Sensitivity analysis of the measured iterations to convergence for the AMG
    preconditioned FGMRES solvers with Vanka and DLSC relaxation. Optimal parameter values, $\eta^{opt}$, are outlined in
    black (every fifth sample is marked).
    }\label{fig:sensitivity_th}
\end{figure}

The optimal damping parameter values, $\eta^{opt}$, used within this
section are shown in~\cref{table:p2p1_taus}. 
To understand this parameter choice, we examine the effect of 
$\eta$ on the convergence rate of the FGMRES solver when using the 4
(3) level AMG
hierarchy preconditioner for 2D (3D) problems with at least 1 million degrees of
freedom on the \isoptwopone{} grid.
The \cMphLO{} preconditioner does not require the selection of $\eta$, as
the FGMRES iterations are invariant to changes in scaling and, in this case,
$\eta$ simply scales the entire preconditioner.
Our results, as shown in~\cref{fig:sensitivity_th}, indicate
that the optimal damping parameters for structured and unstructured 2D
problems using Vanka fall within the range of $0.75-0.86$. The results for
3D problems using Vanka show a steady plateau where the minimum iteration count occurs, with
the optimal $\eta$ parameter for all 3D problems selected as
$1.0$, which falls within this interval.  In contrast, for DLSC, we find $\eta=1.0$ is optimal for 2D problems, but variations are seen for 3D problems. Importantly, across all cases,
underestimating this optimal parameter results in quicker performance
degradation than overestimating it.
\begin{table}[!ht]
\begin{center}
    \begin{tabular}{c c | c c c c}
    \toprule
    & & \multicolumn{4}{c}{Mesh Type} \\
    & & \multicolumn{2}{c}{Structured} & \multicolumn{2}{c}{Unstructured} \\
    \cmidrule(lr){3-4} \cmidrule(lr){5-6}
    Relaxation & Preconditioner & 2D & 3D & 2D & 3D \\
    \midrule
    \multirow{3}{*}{Vanka} & \cMphHO{}   & 0.86 & 1.00 & 0.75 & 1.00 \\
                           & \cMphLO{}   & 1.00 & 1.00 & 1.00 & 1.00 \\
                           & \cMphHLO{}  & 0.75 & 1.00 & 0.80 & 1.00 \\
    \midrule
    \multirow{2}{*}{DLSC}& \cMphHO{}   & 1.00 & 0.60 & 1.00 & 1.00 \\
                           & \cMphHLO{}   & 1.00 & 1.10 & 1.00 & 1.06 \\
    \bottomrule
    \end{tabular}
\end{center}
\caption{Optimal damping parameters $\eta^{opt}=\eta_u^{opt}=\eta_p^{opt}$ for Vanka and DLSC relaxation in Taylor-Hood, \ptwopone{} preconditioners. \cMphLO{} results for DLSC are excluded due to poor performance across all problems.}\label{table:p2p1_taus}
\end{table}

\begin{figure}[!ht]
\centering
\includegraphics[width=\textwidth]{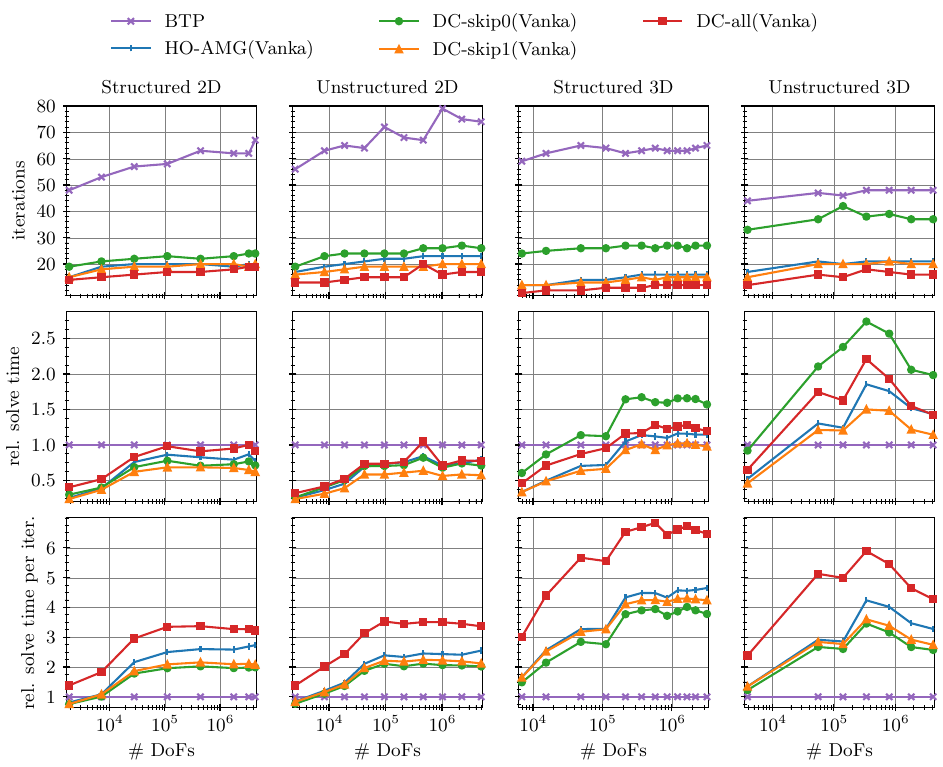}
\captionsetup{skip=1pt}
\caption{\ptwopone{} Vanka relaxation-based monolithic AMG solvers. Iterations to convergence (top),
          Relative time to solution (middle),
          and relative time per iteration (bottom).
          The timings are relative to the \cMBTP{} solver.
          }\label{fig:th_results}
\end{figure}

Using the above-described parameter choices, we examine the convergence of
FGMRES preconditioned with the monolithic AMG-based preconditioners based on Vanka relaxation and compare
them against the \cMBTP{} preconditioner for a range of problems 
in~\cref{fig:th_results}.
For most
problems, the \cMh{} preconditioner (based on the \ptwopone{} AMG hierarchy) exhibits a slight growth in the iteration count as the problems
are refined.
Similar slight growth is seen for the defect-correction (DC) approach based on
the \isoptwopone{} discretization; however,
\cMphHO{}, which effectively replaces the coarse grids of the \cMh{} hierarchy with an AMG
hierarchy based on $K_1$, consistently converges in 1--3 fewer iterations than
\cMh{}. Since the \cMphHO{} coarse grids are sparser than those
of \cMh{}, the relative cost of the solution is slightly cheaper, as can be seen
in the timing rows of~\cref{fig:th_results}. However, the solution times are
similar for \cMh{} and \cMphHO{}, because the relaxation cost at the finest
level dominates the time to solution. For \cMphLO{}, relaxation on the
$K_0$ system is skipped and replaced
by relaxation on the $K_1$ system. The convergence of \cMphLO{}
suffers for all problem types, as seen in the top row
of~\cref{fig:th_results}. Despite having a lower cost per iteration than \cMh{},
the increase in the iteration count makes \cMphLO{} less computationally
effective than \cMphHO{} and \cMh{}. This effect is particularly drastic for 3D
problems, which can be seen in the last two columns of~\cref{fig:th_results}.
Although the solver based
on \cMphHLO{} converges in the fewest number of iterations for most problems, doubling
the number of relaxation sweeps on fine grids (relaxing on both $K_0$ and $K_1$) results
in significantly longer times to solution than for the \cMh{} approach.
However, the relative difference in timing between \cMphHLO{} and \cMh{} appears to
diminish with larger problem sizes.
The convergence and timing results in~\cref{fig:th_results}, especially for
\cMphLO{} and \cMphHO{}, show conclusively that relaxation on the original
problem ($K_0$) should not be skipped when designing a defect-correction type
AMG preconditioner. Additional relaxation on the $K_1$ system may improve
convergence further but at a higher computational cost. This conclusion is
supported by the results from a GMG-based defect-correction approach published
in~\cite{voronin2021low}.

In contrast, the \cMBTP{} preconditioners exhibit slower convergence,
requiring approximately 3--5 times more iterations to reach the convergence tolerance. This
effect is particularly noticeable for 2D problems, where the iteration counts
show significant growth. Comparing the time to convergence, we find that the
monolithic AMG solvers outperform the \cMBTP{}-based solvers by 20--50\% for 2D
problems. For 3D problems, however, \cMBTP{} is generally the fastest to converge, although \cMphHO{} offers similar performance on structured grids.
This difference in performance is partly due to the reduced number of iterations
required for convergence of \cMBTP{} in 3D, coupled with the larger Vanka patch
sizes at each hierarchy level, underscoring the need for further research on
implementation efficiencies.

\begin{figure}[!ht]
\centering
\includegraphics[width=\textwidth]{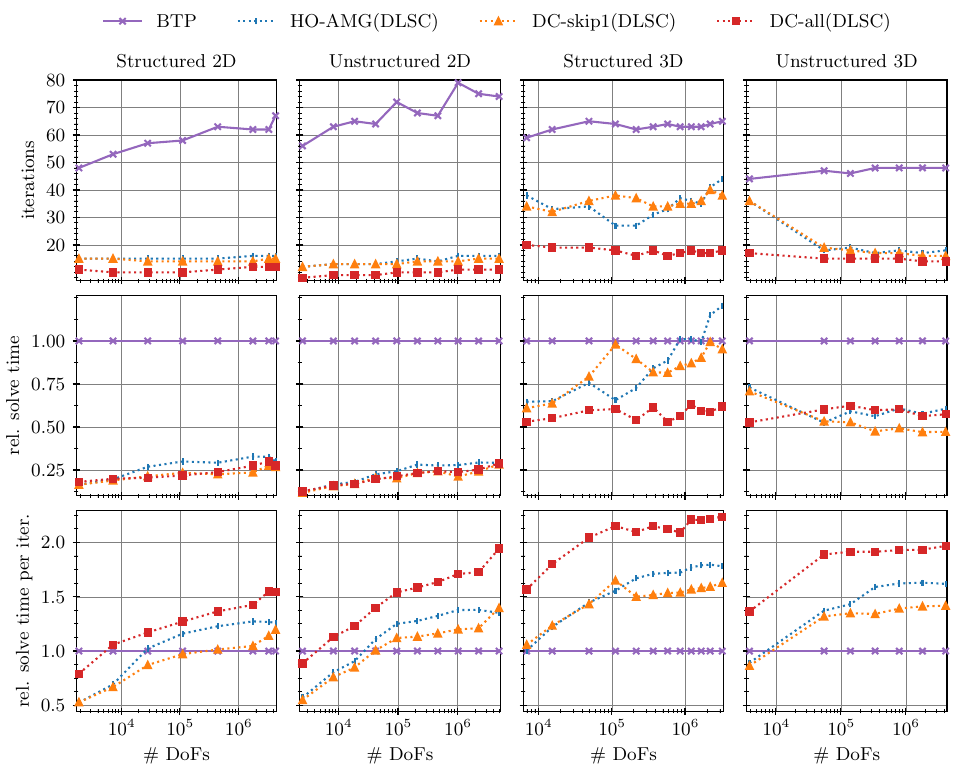}
\captionsetup{skip=1pt}
\caption{\ptwopone{} DLSC relaxation-based monolithic AMG solvers. Iterations to convergence (top),
      Relative time to solution (middle),
      and relative time per iteration (bottom).
      The timings are relative to the \cMBTP{} solver.
}\label{fig:th_lsc_results}
\end{figure}

We next consider the performance of FGMRES preconditioned
with DLSC-based monolithic AMG preconditioners, continuing to compare with
\cMBTP{}. The timing results 
shown in~\cref{fig:th_lsc_results} include results for
\cMBTP{} identical to those in~\cref{fig:th_results}, which helps
evaluate the effectiveness of Vanka versus DLSC-based AMG
solvers.
As seen in the top row of~\cref{fig:th_lsc_results}, the number of iterations
required for the DLSC-based preconditioners follows similar trends to those
observed with Vanka-based AMG solvers. Specifically, the \cMphHLO{}
preconditioners achieve convergence in fewer iterations than both \cMphHO{} and
\cMh{}. The performance of \cMphHO{} and \cMh{} in terms of iterations is
similar across various problem types. A notable difference between Vanka and
DLSC-based AMG preconditioners is the larger increase in iteration count for
3D structured problems in the latter, even after applying additional Chebyshev sweeps as
detailed in~\cref{table:lsd_dgs_config}.

The middle section of~\cref{fig:th_lsc_results} presents an intriguing finding
regarding time-to-solution. For 2D problems, DLSC-based AMG
solvers are 75--80\% faster than \cMBTP{}, converging much more quickly. In 3D scenarios,
this lead diminishes to about 25-50\% for \cMphHO{}, suggesting an average solution speed about 1.3--2 times faster than Vanka-based AMG solvers. Despite similar iteration counts, the defect-correction solver \cMphHO{} is found to be consistently more efficient than \cMh{},
thanks to its sparser coarser operators. The final rows
of~\cref{fig:th_results,fig:th_lsc_results} indicate that, per iteration,
DLSC-based AMG preconditioners are more cost-effective than their Vanka-based
counterparts.

The findings from~\cref{fig:th_results,fig:th_lsc_results} underscore the
importance of including relaxation on the original problem ($K_0$) in
defect-correction type AMG preconditioners. Adding relaxation for $K_1$ improves
convergence but increases computational costs.

\begin{figure}[!ht]
\centering
\includegraphics[width=\textwidth]{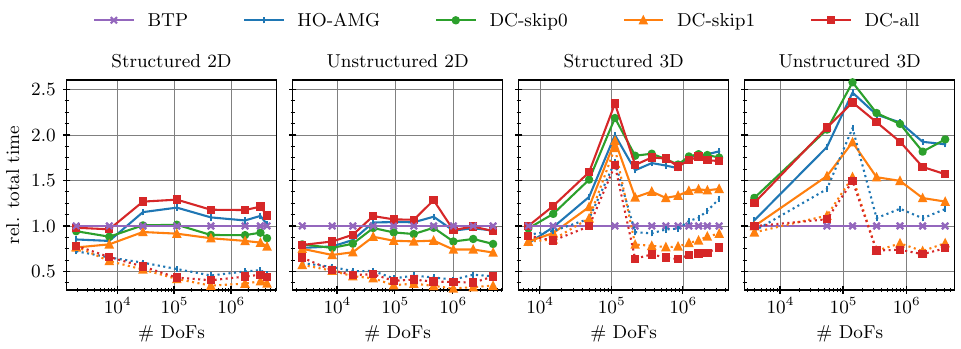}
\captionsetup{skip=1pt}
\caption{Relative total time to solution, including multigrid setup and solve phases, for \ptwopone{} problems. The solid HO and DC lines correspond to solvers utilizing Vanka relaxation and dashed lines correspond to LSC-DGS relaxation. The timings are relative to the \cMBTP{} solver.
}\label{fig:th_total_timings}
\end{figure}

Turning to the total time to solution, which includes both multigrid setup and solve phases, \Cref{fig:th_total_timings} summarizes the performance of all solver types described thus far.  The plots clearly show that the LSC-DGS-based multigrid methods (dashed lines) are consistently faster than Vanka-based solvers (solid lines), a difference attributed to their faster setup and solve phases. For 2D problems, most Vanka-based solvers and all LSC-DGS-based solvers outperform \cMBTP{}. However, in the 3D case where the cost of Vanka assembly and application rises drastically, none of the Vanka-based solvers can outperform \cMBTP{} for large problem sizes. In contrast, the LSC-DGS-based solvers outperform \cMBTP{} by 10-30\% for most large 3D problem sizes. %
Overall, for all problems considered, DLSC-based AMG preconditioners offer significant advantages in speed and cost-effectiveness over Vanka-based methods.

\subsection{Preconditioning~\texorpdfstring{\ptwoponedisc}{P2/P1disc}
Systems}\label{sec:p2p1disc_results}

\begin{figure}[!ht]
    \centering
    \includegraphics[width=\textwidth]{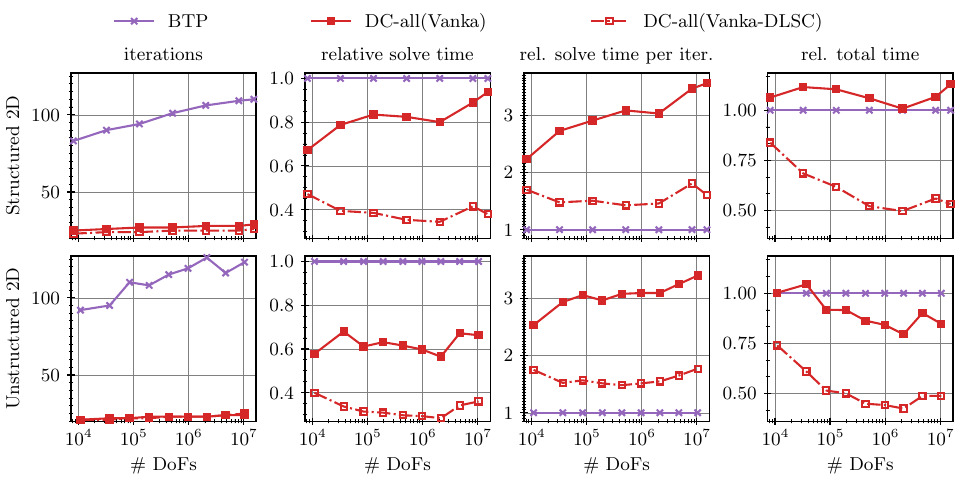}%
    \caption{\ptwoponedisc{} solvers: Iterations to convergence (leftmost), and relative metrics compared to \cMBTP{} solver - time to solution (middle left), time per iteration (middle right), and total time (rightmost), including multigrid setup and solve phase times.
    }\label{fig:sv_results}
\end{figure}

In this section, we evaluate the effectiveness of the defect-correction
preconditioner for solving systems that have been discretized using
\ptwoponedisc{} elements on 2D structured and unstructured meshes.  Preliminary
numerical experiments (not reported here) showed that relaxation on all levels
of the AMG hierarchy is needed to generate an efficient solver.
The convergence data in~\cref{fig:sv_results} uses the
optimal damping parameter values tabulated in~\cref{table:p2p1disc_params}.
\begin{table}[!ht]
\begin{center}
\begin{tabular}{c | c c c}
\toprule
           &             & \multicolumn{2}{c}{Parameters}\\
           \cmidrule(lr){3-4}
Relaxation & Mesh Type   &$\omega_0^{opt}$  &$\eta_p^{opt}$   \\
\midrule
\multirow{2}{*}{Vanka($\ell\geq 0$)}   
						 & Structured    &0.78       &2.68       \\ 
                         & Unstructured  &0.58       &4.00       \\ 
\midrule
\multirow{2}{*}{Vanka($\ell=0$),DLSC($\ell>0$)}  
						 & Structured    &0.87       &2.96       \\ 
                         & Unstructured  &0.58       &3.80       \\ 
\bottomrule
\end{tabular}
\end{center}
\caption{Optimal parameters for Scott-Vogelius, \ptwoponedisc{}, \cMphHLO{} preconditioners: overweighting $\eta_p^{opt}$, and damping $\omega_0^{opt}$.} \label{table:p2p1disc_params}
\end{table}

We study the sensitivity of the FGMRES convergence rate to
the parameter choices in \cMphHLO{} by performing a parameter search on 3-level AMG
hierarchies with at least a quarter million degrees of freedom on the finest
\isoptwopone{} grid for both structured and unstructured problems.
We explore two \cMphHLO{} defect-correction strategies: the first applying
Vanka relaxation across all levels except the coarsest, while the second, termed
``Vanka-DLSC,'' utilizes Vanka relaxation at the finest level ($\ell=0$) and
DLSC on coarser levels ($\ell > 0$). Both strategies adhere to the relaxation
sweep schedules previously discussed.

\begin{figure}[!ht]
    \centering
    \begin{subfigure}{\textwidth}
    	\centering
    	\includegraphics[width=0.8\textwidth]{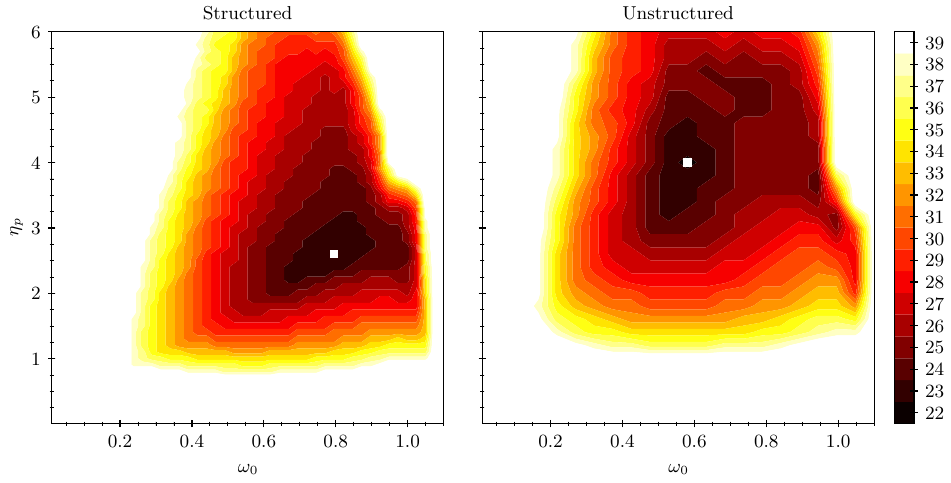}
    	\captionsetup{skip=1pt}
        \caption{AMG preconditioners using Vanka relaxation.}\label{fig:sensitivity_sv_vanka}
    \end{subfigure}
    \begin{subfigure}{\textwidth}
        \centering
        \includegraphics[width=0.8\textwidth]{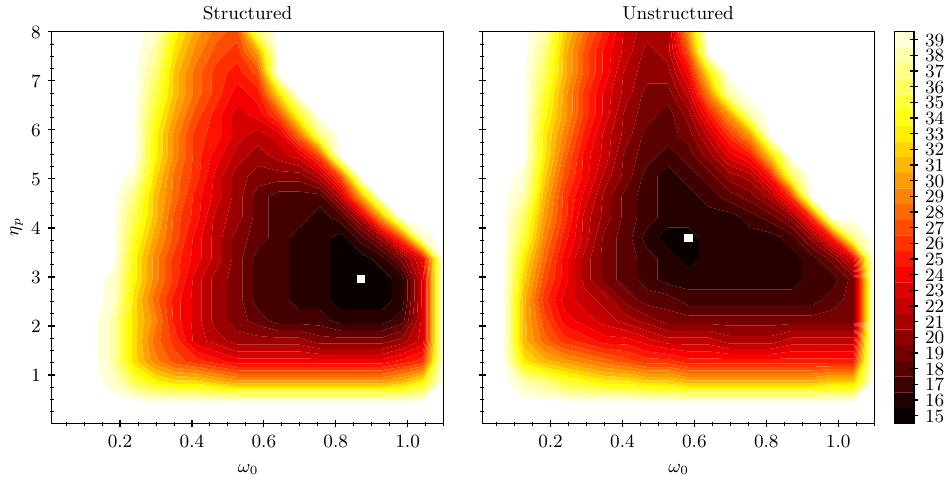}\
        \captionsetup{skip=1pt}
        \caption{AMG preconditioners using Vanka on $\ell=0$ and DLSC on $\ell>0$ multigrid levels.}\label{fig:sensitivity_sv_lsc_dgs}
    \end{subfigure}
    \captionsetup{skip=1pt}
	\caption{Sensitivity analysis of iterations to convergence for \cMphHLO{} preconditioned FGMRES solvers. Optimal parameters are highlighted with white squares.}    
    \label{fig:sensitivity_sv}
\end{figure}

We found that for both relaxation types
$\eta_u \approx 1$ produced the best convergence results, regardless of the
values of $\omega_0$ and $\eta_p$, and so we fix $\eta_u^{opt}=1$.
To determine the other parameters, we perform a 2D parameter scan.
The two-dimensional heat maps
in~\cref{fig:sensitivity_sv} show the iteration counts for different choices
of $\omega_0$ and $\eta_p$. The results indicate that the iteration counts
are relatively insensitive to the choice of $\omega_0$ as long as they are
between $0.5$ and $0.9$, and a similar trend is observed for the overweighting
parameter, $\eta_p$, where values between $1.8$ and $5$ result in
the best convergence. The optimal convergence regions for structured and
unstructured problems are relatively close, but unstructured problems tend to
perform best with smaller $\omega_0$ and larger $\eta_p$ values than
structured problems. The optimal parameter choices are summarized
in~\cref{table:p2p1disc_params}.

Using the parameter choices described above, we now examine the
convergence of these solvers across a range of problems and problem sizes.
The iteration counts in~\Cref{fig:sv_results} show that both \cMphHLOVanka{} and \cMphHLOLSC{}
preconditioners converge, on average, in 27 iterations for structured domains and
23 iterations for unstructured domains, with the number of iterations increasing by
one roughly with each order of magnitude increase in problem size.
It is important to recall that the \cMh{} method does not converge on a
significant number of the \ptwoponedisc{} test problems and, so, it is
revealing that \cMphHLO{} converges on all test problems.

In contrast, the \cMBTP{} preconditioners require around 4--5 times more
iterations to reach convergence, and the iteration counts exhibit more
significant growth. With each refinement of the problem (increasing the problem
size by $2\times$), an additional 2--4 iterations are needed in \cMBTP{}.
Despite the substantial difference in iteration counts between \cMBTP{} and
\cMphHLOVanka{}, the overall solve time of \cMBTP{} is only 10--30\% slower for
structured problems and 30--40\% slower for unstructured problems. 
However, compared to \cMphHLOLSC{}, this relative solve time of
\cMBTP{} is 50--65\% slower for structured problems and 60--70\% for unstructured
grids. This 1.5--2$\times$ improvement in time to solution over \cMphHLOVanka{} suggests
better solve-phase cost-effectiveness in using \cMphHLOLSC{} for \ptwoponedisc{}
discretizations, which is supported by the relative solve phase timings per iteration.
The last column of~\cref{fig:sv_results} examines the relative total times to solution, including multigrid setup and solve phase timings.
These plots are similar to the relative solve phase timings, but the values are shifted upwards relative to the \cMBTP{} line. This upward shift indicates that the relative setup phase for the DC solvers is longer than for \cMBTP{}. As a result, \cMphHLOVanka{} outperforms \cMBTP{} only in the unstructured case, while \cMphHLOLSC{} still outperforms \cMBTP{} in all cases.
It is important to note that the timing results associated with the setup and solve phases are highly dependent on the implementation of relaxation setup and solution phase kernels. These timings are expected to vary depending on the specific implementation and the compute architecture.

Our results demonstrate the effectiveness of the \cMphHLO{} preconditioner for
solving systems discretized using the \ptwoponedisc{} discretizations. These
findings validate the capabilities of our approach and offer valuable guidance for
selecting optimal parameter values when utilizing the solver. Importantly, the
comparison with the \cMBTP{} method reveals the superiority of our approach in
terms of iteration count and runtime, with the \cMphHLOLSC{} showing promising
performance improvements over Vanka-based monolithic AMG preconditioners. 

\section{Application Problems}\label{sec:application_problems}
\begin{figure}[!ht]
    \centering
    \begin{subfigure}{0.48\textwidth}
        \centering
        \includegraphics[width=\textwidth]{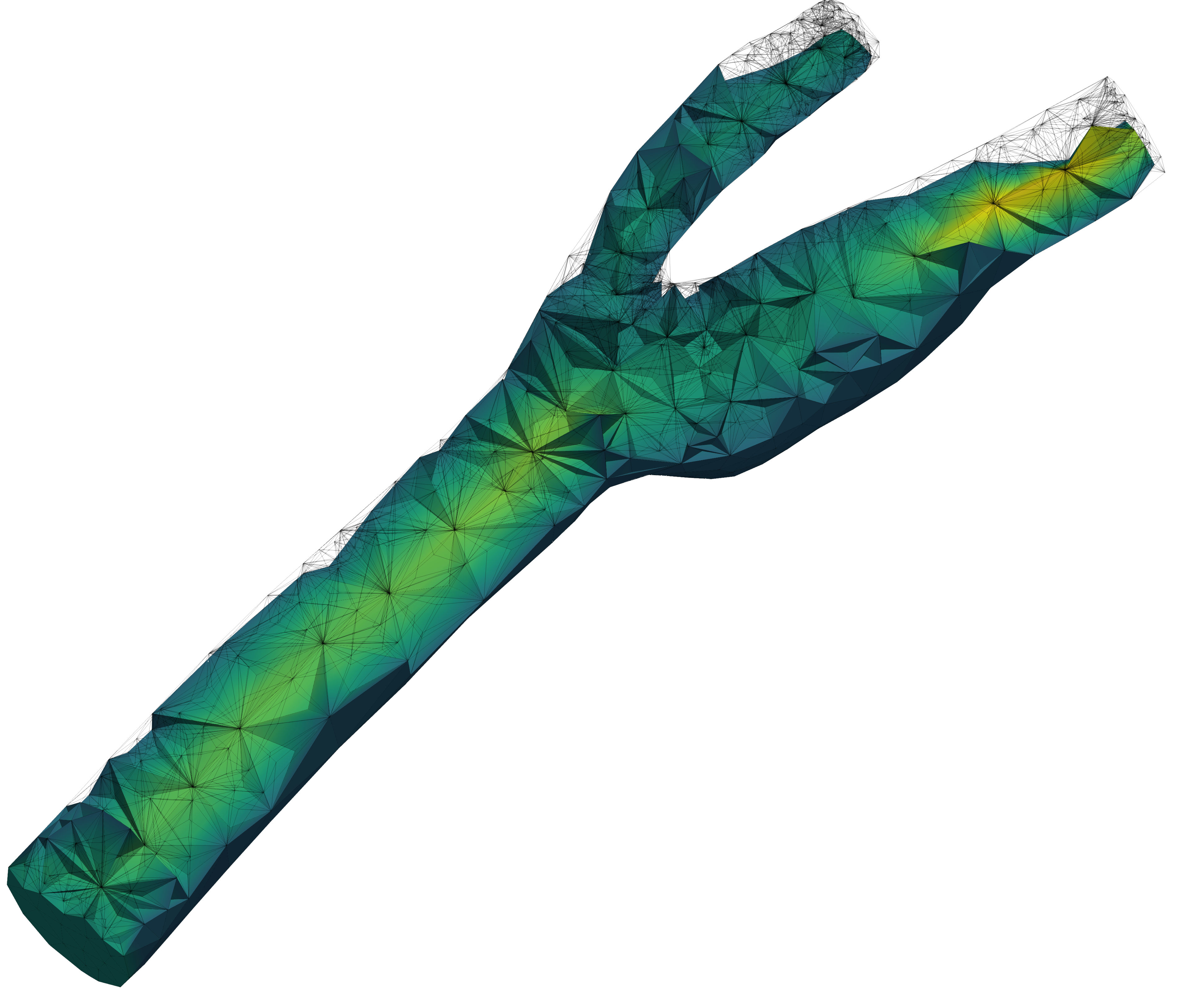}
        \caption{3D Artery: \ptwopone{}}\label{fig:artery}
    \end{subfigure}
    \hfill
    \begin{subfigure}{0.4\textwidth}
        \centering
        \includegraphics[width=\textwidth]{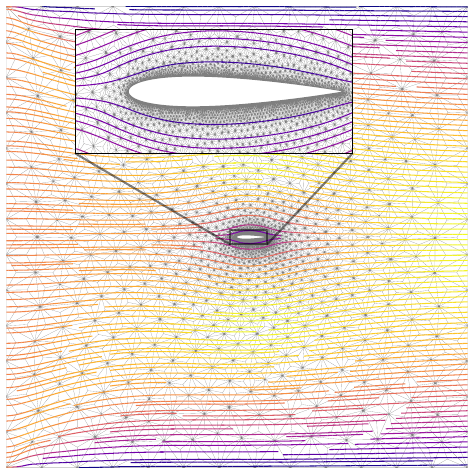}
        \caption{2D Airfoil: \ptwoponedisc{}}\label{fig:airfoil}
    \end{subfigure}
    \hfil
    \captionsetup{skip=1pt}
    \caption{Test problems for Section~\ref{sec:application_problems}}\label{fig:applications_problems}
\end{figure}

We conclude by examining the preconditioned FGMRES
algorithms on two realistic applications pictured in~\cref{fig:applications_problems}:
\begin{description}
\item[3D Artery: \ptwopone{}] \hfill \\
    Models fluid flow inside an artery~\cite{geuzaine2009gmsh},
    discretized with unstructured tetrahedral elements.
    We used an STL file obtained from a luminal casting of a carotid artery
    bifurcation\footnote{Carotid artery bifurcation CAD credit
    \url{https://grabcad.com/library/carotid-bifurcation}}
    to generate a series of meshes of increasing
    sizes~\cite{geuzaine2009gmsh}. The inflow
    BC is imposed normal to the bottom left cross-section of
    the artery (see~\cref{fig:artery}), and the outflow BC
    (Neumann) are imposed at the opposite two cross-sections,
    with zero Dirichlet velocities on the 
    other boundary faces of the domain.
\item[2D Airfoil: \ptwoponedisc{}] \hfill \\
    This problem models flow around an airfoil~\cite{geuzaine2009gmsh} in
    the center of the domain, with an adaptively refined 2D unstructured
    mesh of triangular elements\footnote{Airfoil \textit{geo} file is
    available in Gmsh's GitHub~\href{https://github.com/jeromerobert/gmsh/blob/b46d5e673c48f5108b373014aabf9d69e970a80f/benchmarks/3d/naca0012.geo}{repository}}. A parabolic inflow BC is imposed along the
    left edge of the domain, while a homogeneous Neumann BC is imposed on the outflow boundary at the right edge,
    with zero Dirichlet BCs on the velocity field
    on all other boundary edges of the domain (including the inner airfoil).
\end{description}
In all tests, we adhere to the convergence criteria and data collection methods
outlined in~\cref{sec:results}. The Vanka
relaxation-based AMG results reported in this section use
the same solver parameters described
in~\cref{table:p2p1_taus,table:p2p1disc_params}. However, the DLSC
relaxation-based preconditioners required enhanced approximation stages for the
continuity equation and pressure correction~---~see
lines \ref{line:continuity} and \ref{line:pressure_correction} of~\cref{alg:DLSC}. For both 2D Scott-Vogelius and 3D Taylor-Hood
problems, we substitute Chebyshev polynomial relaxation in these stages with
SA-AMG V-cycle using weighted Jacobi relaxation. This does not significantly
affect the optimal damping parameters $(\omega_0^{opt},\eta_p^{opt})$,
demonstrating the parameter-free nature of the DLSC preconditioners.
Additionally, for the 3D problem, increasing the V-cycle's pre- and
post-relaxation sweeps from two to three iterations results in improved
convergence rates.

\begin{figure}[!ht]
\centering
\includegraphics[width=\textwidth]{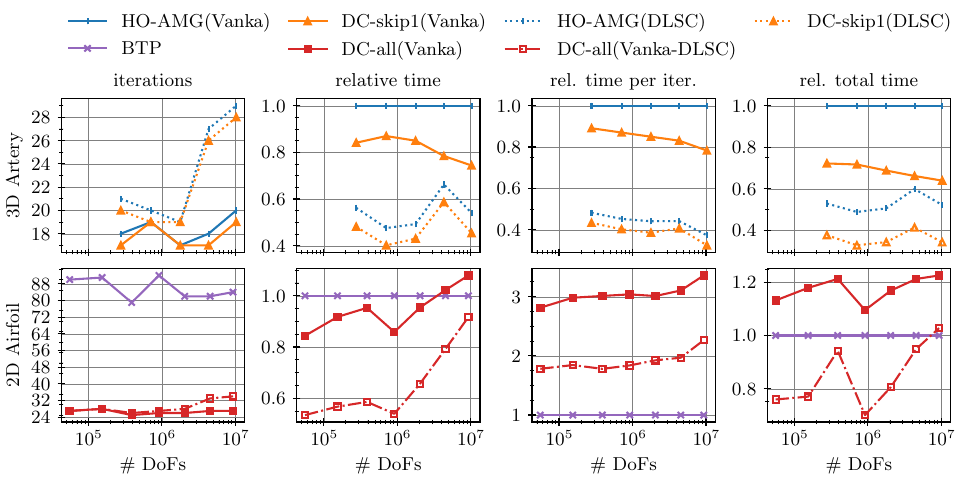}
\caption{Top row: 3D \ptwopone{} solvers; Bottom row: 2D \ptwoponedisc{} solvers. 
Columns show iterations to convergence (leftmost) and relative metrics compared to a reference solver - time to solution (middle left), time per iteration (middle right), and total time (rightmost), including multigrid setup and solve phase times.
Timings are relative to \cMh{} for 3D and \cMBTP{} for 2D.
}\label{fig:hero_conv}
\end{figure}

\Cref{fig:hero_conv} showcases the iteration counts for multigrid preconditioned FGMRES, accompanied by the relative time per solution and per iteration for each
solver type. The top row of~\cref{fig:hero_conv} is associated with the 3D
artery problem, discretized with the \ptwopone{} element pair, whereas the
bottom row corresponds to the 2D airfoil problem, discretized with
\ptwoponedisc{} elements.
For the 3D problem, we compare \cMphHO{} timings against \cMh{} instead of 
comparing with \cMBTP{}.
For elongated domains, the pressure mass
matrix is ineffective as a preconditioner for the pressure Schur
complement~\cite{AJanka_2008a,dobrowolski2003lbb}, leading to poor convergence
of \cMBTP{} in the 3D artery problem, where \cMBTP{}
requires 900--1500 iterations, taking 20--30 times longer than
\cMphHO{}. Hence, we instead compare \cMphHO{} to \cMh{}.

For the 3D problem discretized with \ptwopone{}, \cMphHO{} solvers converge in fewer or equivalent iterations compared to \cMh{} with the same relaxation type. 
Thus, the relative efficiency of \cMphHO{} is attributed to both reduced
iteration numbers and faster individual iterations. The \isoptwopone{} coarse
grids in \cMphHO{} are sparser than \ptwopone{} grids, leading to more compact
Vanka patches and decreased computational demand. Consequently, \cMphHO{}
solvers can handle larger problems with similar hardware resources. 
The performance gap is more noticeable in Vanka-based solvers due to their higher
storage and computational requirements. \cMphHOVanka{} and \cMphHODLSC{} solvers
are approximately 13--25\% and 2--10\% faster, respectively, than their
\cMh{} counterparts. 
Despite taking more iterations to converge, DLSC-based preconditioners are faster than Vanka-based
preconditioners due to the much cheaper cost per iteration, even with higher iteration counts. However, as noted in~\cref{rmrk:dlsc}, DLSC preconditioners may not be suited for all complex geometries. 
The relative total time column for the 3D \ptwopone{} problem follows similar
trends to the relative time to solution, but with all lines shifted by 10\% lower relative to \cMh{}. This indicates that even the setup for the \cMphHO{} is faster due to the sparser coarse grid. 

In the analysis of the two-dimensional \ptwoponedisc{} airfoil problem, we
conducted a comparative study between low-order \cMphHLO{} preconditioners and
the \cMBTP{} approach. The preconditioners \cMphHLOVanka{} and \cMphHLOLSC{}
demonstrated convergence in significantly fewer iterations than
\cMBTP{}. In particular, for smaller-scale problems, \cMphHLOVanka{} achieved solutions 5--15\% more quickly than \cMBTP{}. However, as problem size increased, this advantage diminished, leading to a performance reduction of up to 8\%. In contrast, \cMphHLOLSC{} consistently outperformed \cMBTP{}, with solution times 40-45\% faster for smaller problems, but with decreased efficiency for larger problem sizes. As the problem size increases and multigrid hierarchies grow deeper, \cMphHLOLSC{}
becomes less competitive, and its convergence rates start diverging from those
of \cMphHLOVanka{}.
A consistent trend emerges when comparing these results to those obtained
from two-dimensional unstructured meshes (see Figure~\ref{fig:sv_results}). For
smaller problem sizes, the cost per iteration of both \cMphHLOLSC{}
preconditioners, compared to \cMBTP{}, remains relatively stable. However,
a noticeable increase in cost is observed as problem sizes increase. This
is attributed to denser $BM_u^{-1}B^T$ operators on coarser grids and a reduced computational efficiency of SA-AMG in preconditioning these operators.
The relative total time column displays similar behavior as observed in~\cref{sec:p2p1disc_results}. 
These plots are similar to the relative solve phase timings, but the values are shifted upwards relative to the \cMBTP{} line due to higher setup costs. As a result, \cMphHLOVanka{} is always slower than \cMBTP{}, while \cMphHLOLSC{} outperforms \cMBTP{} in all but the largest case.

\begin{remark}[DLSC relaxation-based AMG Preconditioners]\label{rmrk:dlsc}
In experiments not reported here, we observe that AMG solvers with DLSC
relaxation often do not converge (within 200 iterations) for certain arterial
simulations, a contrast to Vanka relaxation-based solvers. Attempts to enhance the DLSC solver,
such as increasing the number of sweeps in DLSC stages or boosting
the count of pre and post-relaxation iterations, results in negligible improvements.
The difference in performance is largely attributed to mesh quality issues,
notably the presence of poorly shaped elements. To address this, we incorporated
10 iterations of Laplacian mesh smoothing during mesh generation in the present
study, which improves the results. Despite this intervention, DLSC
solvers still demonstrate more variability in iteration counts than Vanka solvers across problem types.  
\end{remark}

\begin{remark}[Velocity Mass Matrix Scaling in DLSC]\label{rmrk:mass_mat_scaling}
The scaling of the velocity mass matrix, as outlined in~\cref{alg:DLSC}, led to reduced convergence rates for the Taylor-Hood discretized artery problem solvers, as depicted in the top row of~\cref{fig:hero_conv}. This issue predominantly arises in the presence of challenging anisotropic meshes. In an attempt to address this, our experiments with direct inversion of the velocity mass matrix for smaller-scale, two-level hierarchies did not significantly enhance convergence rates. Remarkably, for the artery mesh, the approximation of $M_u=I$ consistently resulted in the lowest interaction counts and the fastest solution times. These findings underscore the need for a more comprehensive investigation into the implications of mass matrix scaling within the DLSC algorithm on complex meshes. This remains an important area for future research.
\end{remark}

\section{Conclusion}\label{sec:conclusion}

We present a novel algorithm for constructing monolithic
algebraic multigrid (AMG) preconditioners for the Stokes problem using a
defect-correction approach based on the lower-order \isoptwopone{}
discretization. By splitting the Stokes operator into
block-matrix format, we simplify the assembly of algebraic interpolation
operators for each variable type, thus boosting efficiency. Our approach employs
the pressure stiffness matrix or an algebraic equivalent, enhancing solver robustness.
This methodology can precondition both the
\ptwopone{} and \ptwoponedisc{} discretizations. When applied directly to the \ptwopone{}
discretization, this framework also results in a robust solver, albeit at a slightly higher 
iteration cost than low-order preconditioners.  Both methods
consistently demonstrate a relatively steady iteration count across various
Stokes problems on structured and unstructured meshes in 2D and 3D. 
Notably, the monolithic AMG preconditioners often match or even surpass the
performance of inexact block triangular preconditioners across a range of
problems. The robustness of the monolithic AMG approach is particularly
pronounced for the elongated-domain problem, which is problematic
for the inexact block triangular preconditioner.
An important advance in this work is the novel block factorization-based
solution for Vanka patch systems, significantly reducing storage and
computational expenses up to $2\times$ in 2D and $3\times$ in 3D. We also introduce
a Chebyshev-based adaptation of the LSC-DGS~\cite{wang2013multigrid}
relaxation aimed at enhancing parallelism.
This relaxation approach proves to be highly effective and efficient in well-meshed 
domains but less resilient in complex ones, unlike Vanka relaxation.
Our AMG preconditioners with Vanka relaxation use FGMRES to automatically damp relaxation
for broader applicability without
intricate tuning. DLSC relaxation-based AMG preconditioners, while not needing
damping parameters do require appropriate solver and sweep count choices for
each relaxation stage, which is subject to problem-dependent tuning. 
For the Scott-Vogelius discretization, our findings show that optimal
parameters associated with the defect-correction framework remain stable under
minor perturbations and can be efficiently determined through a coarse-parameter
scan.

\section{Future Work}\label{sec:future_work}

We aim to extend the \isoptwopone{}-based defect-correction preconditioners to
higher-order Taylor-Hood and Scott-Vogelius discretizations,
by employing a $p$-multigrid approach.
Considering the inherent challenges of such discretizations, these
preconditioners are expected to be more advantageous than directly applying AMG
to high-order systems.
In the context of higher-order solvers, efficient Vanka implementations on fine
levels are set to become an increasingly vital component of the described
defect-correction approach. With advancements in supercomputing, these
higher-order Vanka methods are expected to be highly suitable for many-core
systems, thanks to their expansive stencils and uniform structure. In contrast,
low-order coarse-level Vanka patches, characterized by irregular dimensions and sparsity patterns,
pose significant challenges in designing efficient data structures. This
complexity at the coarse level underscores the utility of DLSC relaxation. DLSC
relaxation, known for its efficiency and reduced memory demands, offers a
practical approach to avoid the irregularities and computational hurdles
associated with low-order coarse-level Vanka patches.
However, it should be noted that DLSC relaxation also presents specific challenges, particularly with complex unstructured meshes, which can adversely affect convergence rates as discussed in~\cref{rmrk:mass_mat_scaling}.
Additionally, we plan to strengthen the theoretical foundation of our methods,
aiming to derive convergence guarantees for the proposed monolithic AMG
preconditioner. A key step will be ensuring the stability of coarsened velocity
and pressure DoFs, essential for maintaining the discrete inf-sup condition and
reinforcing our methodologies' reliability.

\section*{Acknowledgments}

The work of SPM was partially supported by an NSERC Discovery Grant. RT was
supported by the U.S.~Department of Energy, Office of Science, Office of
Advanced Scientific Computing Research, Applied Mathematics program. Sandia
National Laboratories is a multi-mission laboratory managed and operated by
National Technology \& Engineering Solutions of Sandia, LLC (NTESS), a wholly owned subsidiary of Honeywell International Inc., for the U.S. Department of Energy’s National Nuclear Security Administration (DOE/NNSA) under contract DE-NA0003525. This written work is authored by an employee of NTESS. The employee, not NTESS, owns the right, title and interest in and to the written work and is responsible for its contents. Any subjective views or opinions that might be expressed in the written work do not necessarily represent the views of the U.S. Government. The publisher acknowledges that the U.S. Government retains a non-exclusive, paid-up, irrevocable, world-wide license to publish or reproduce the published form of this written work or allow others to do so, for U.S. Government purposes. The DOE will provide public access to results of federally sponsored research in accordance with the DOE Public Access Plan.

\bibliographystyle{siamplain}
\bibliography{references}
\end{document}